\documentclass[epic,11pt]{article}
\usepackage{amsfonts,amsmath,amscd,amssymb,amsthm}
\usepackage[alphabetic,initials]{amsrefs}
\usepackage[all]{xy}
\usepackage{hyperref} 
\usepackage{mathrsfs}
\input xy
\usepackage{fullpage}
\usepackage{amssymb}
 
\def\convf{\hbox{\space \raise-2mm\hbox{$\textstyle      \bigotimes \atop \scriptstyle \omega$} \space}}
\def\0{{\bar 0}}
\def\1{{\bar 1}}

\def\Z{{\mathbb Z}}

\def\B{{\mathbb B}}

\def\N{{\mathbb N}}

\def\soc{{\operatorname{soc}}}

\def\Ind{{\operatorname{Ind}}}

\def\pwg1{{\operatorname{PWG}}}
\def\pwg{{\operatorname{pwg}}}

\def\span{{\operatorname{span}}}

\def\Hom {{\operatorname{Hom}}}
\def\Ker {{\operatorname{Ker}\;}}

\newcommand{\ttk}{\mathtt{k}}
\newcommand{\tte}{\mathtt{e}}

\newcommand{\gL}{\Lambda}

\newcommand{\itema}{\item[{{\rm(a)}}]}
\newcommand{\itemb}{\item[{{\rm(b)}}]}
\newcommand{\itemc}{\item[{{\rm(c)}}]}
\newcommand{\itemd}{\item[{{\rm(d)}}]}
\newcommand{\iteme}{\item[{{\rm(e)}}]}
\newcommand{\itemf}{\item[{{\rm(f)}}]}

\newcommand{\noi}{\noindent}
\newcommand{\ga}{\alpha}
\newcommand{\gb}{\beta}
\newcommand{\gc}{\gamma}

\newcommand{\Gd}{\Delta}

\newcommand{\gd}{\delta}

\newcommand{\gs}{\sigma}

\newcommand{\gO}{\Omega}

\newcommand{\gt}{\tau}
\newcommand{\gz}{\zeta}

\newcommand{\gl}{\lambda}

\newcommand{\gr}{\rho}

\newcommand{\gep}{\epsilon}
\newcommand{\gth}{\theta}
\newcommand{\op}{\oplus}
\newcommand{\gve}{\varepsilon}

\def\Im{{\operatorname{Im}\;}}

\newcommand{\ot}{\otimes}

\newcommand{\fg}{\mathfrak{g}}
\newcommand{\fsl}{\mathfrak{sl}}\newcommand{\fpsl}{\mathfrak{psl}}\newcommand{\osp}{\mathfrak{osp}}

\newcommand{\fr}{\mathfrak{r}}

\newcommand{\ft}{\mathfrak{t}}
\newcommand{\fh}{\mathfrak{h}}
\newcommand{\fz}{\mathfrak{z}}
\newcommand{\fb}{\mathfrak{b}}
\newcommand{\fm}{\mathfrak{m}}
\newcommand{\fn}{\mathfrak{n}}

\newcommand{\fk}{\mathfrak{k}}
\newcommand{\fp}{\mathfrak{p}}
\newcommand{\fq}{\mathfrak{q}}
\newcommand{\fl}{\mathfrak{l}}

\newcommand{\ff}{\footnote}
\newfont{\eufm}{eufm10 scaled\magstep1}

 \newcommand{\ti}{\times}

\newcommand{\cO}{\mathcal{O}}

\newcommand{\cE}{\mathcal{E}}

\newcommand{\cH}{\mathcal{H}}

\newcommand{\ey}{\end{eqnarray}}
\newcommand{\by}{\begin{eqnarray}}
\newcommand{\nn}{\nonumber}

\newcommand{\bco}{\begin{conjecture}}
\newcommand{\ba}{\begin{alg}}
\newcommand{\ea}{\end{alg}}
\newcommand{\eco}{\end{conjecture}}
\newcommand{\bpf}{\begin{proof}}
\newcommand{\epf}{\end{proof}}
\newcommand{\bt}{\begin{theorem}}

\newcommand{\et}{\end{theorem}}

\newcommand{\br}{\begin{rem}}
\newcommand{\er}{\end{rem}}
\newcommand{\brs}{\begin{rems}}
\newcommand{\ers}{\end{rems}}
\newcommand{\bi}{\begin{itemize}}
\newcommand{\ei}{\end{itemize}}
\newcommand{\bl}{\begin{lemma}}
\newcommand{\bsul}{\begin{sublemma}}
\newcommand{\esul}{\end{sublemma}}
\newcommand{\bp}{\begin{proposition}}
\newcommand{\be}{\begin{equation}}
\newcommand{\bc}{\begin{corollary}}
\newcommand{\bexs}{\begin{examples}}
\newcommand{\eexs}{\end{examples}}
\newcommand{\bexa}{\begin{example}}
\newcommand{\eexa}{\end{example}}
\newcommand{\bex}{\begin{exercise}}
\newcommand{\eex}{\end{exercise}}
\newcommand{\btab}{\begin{tab}}
\newcommand{\etab}{\end{tab}}
\newcommand{\bg}{\begin{fig}}
\newcommand{\eg}{\end{fig}}
\newcommand{\el}{\end{lemma}}
\newcommand{\ep}{\end{proposition}}
\newcommand{\ee}{\end{equation}}
\newcommand{\ec}{\end{corollary}}
\newcommand{\Bc}{\begin{center}}
\newcommand{\Ec}{\end{center}}
\newcommand{\bh}{\begin{hyp}}
\newcommand{\eh}{\end{hyp}}
\newcommand{\bhs}{\begin{hyps}}
\newcommand{\ehs}{\end{hyps}}
\newcommand{\bd}{\begin{dfn}}
\newcommand{\ed}{\end{dfn}}

\begin{document}
\title{Table of Contents}

\newtheorem{thm}{Theorem}[section]
\newtheorem{hyp}[thm]{Hypothesis}
 \newtheorem{hyps}[thm]{Hypotheses}
  \newtheorem{rems}[thm]{Remarks}
\newtheorem{conjecture}[thm]{Conjecture}
\newtheorem{theorem}[thm]{Theorem}
\newtheorem{theorem a}[thm]{Theorem A}
\newtheorem{example}[thm]{Example}
\newtheorem{examples}[thm]{Examples}
\newtheorem{corollary}[thm]{Corollary}
\newtheorem{rem}[thm]{Remark}
\newtheorem{lemma}[thm]{Lemma}
\newtheorem{sublemma}[thm]{Sublemma}
\newtheorem{cor}[thm]{Corollary}
\newtheorem{proposition}[thm]{Proposition}
\newtheorem{exs}[thm]{Examples}
\newtheorem{ex}[thm]{Example}
\newtheorem{exercise}[thm]{Exercise}
\numberwithin{equation}{section}%
\setcounter{part}{0}
\newcommand{\drar}{\rightarrow}
\newcommand{\lra}{\longrightarrow}\newcommand{\lda}{\downarrow}
\newcommand{\rra}{\longleftarrow}
\newcommand{\dra}{\Rightarrow}
\newcommand{\dla}{\Leftarrow}
\newcommand{\lfa}{\leftrightarrow}
\newtheorem{Thm}{Main Theorem}
\newtheorem*{thm*}{Theorem}
\newtheorem{lem}[thm]{Lemma}
\newtheorem{fig}[thm]{Figure}
\newtheorem*{lem*}{Lemma}
\newtheorem*{prop*}{Proposition}
\newtheorem*{cor*}{Corollary}
\newtheorem{dfn}[thm]{Definition}
\newtheorem*{defn*}{Definition}
\newtheorem{notadefn}[thm]{Notation and Definition}
\newtheorem*{notadefn*}{Notation and Definition}
\newtheorem{nota}[thm]{Notation}
\newtheorem*{nota*}{Notation}
\newtheorem{note}[thm]{Remark}
\newtheorem*{note*}{Remark}
\newtheorem*{notes*}{Remarks}
\newtheorem{hypo}[thm]{Hypothesis}
\newtheorem*{ex*}{Example}
\newtheorem{prob}[thm]{Problems}
\newtheorem{conj}[thm]{Conjecture}
\title{How to Construct the 
Lattice of Submodules of a Multiplicity free Module from Partial Information.}
\author{Ian M. Musson \\Department of Mathematical Sciences\\
University of Wisconsin-Milwaukee\\ email: {\tt
musson@uwm.edu}}
\maketitle
\begin{abstract} 
In general it is a difficult problem to construct the lattice of submodules 
$L(M)$ of a given module $M$. In \cite{St} R. P. Stanley  outlined  a method for constucting a distributive lattice from a knowledge of its join irreducibles. However it is not an  easy task to identify all join irreducible submodules of a given module. In the case 
of a multiplicity free module $M$ we present a modifiiction of Stanley's method based on the  composition factors of $M$.  As input we require a set of submodules 
$A_1,\ldots , A_n$ whose submodule lattices are known and which contain all composition factors of $M$.  From this we can reconstruct $L(M)$. We illustrate the process for a family of Verma modules $M(\gl_n)$, with $n $ a positive integer, for the Lie superalgebra $\osp(3,2)$.  We show that for $n\ge 2$, $L(M(\gl_n))$ is isomorphic to the (extended) free distributive lattice of rank 3. 
\end{abstract}
\noi \ref{uv.1}. {Introduction.} \\
 \ref{osp}. {Lie Superalgebras.}\\
\ref{elc}. {A Family of  Elusive Cases.}\\
\ref{osp1}. {The Cases $n\ge 2$.}\\
 \ref{uv2}. {Application: Building the Lattice of Submodules of  the Verma  Modules $M(\gl_n)$.} 
\section{Introduction.} \label{uv.1} 
The extended free distributive lattice of rank 3  is a well-known, but quite complicated lattice.  The modules $M(\lambda_n)$ from the abstract have 
 8 composition factors, each with multiplicity one.  For $n\ge 2$, $M(\lambda_n)$ has  20 submodules.  On the other hand  $M(\lambda_1)$ has 23 submodules. 
Studies of $L(N)$ for a module $N$ that arises naturally in representation theory, seem to be quite rare. It is probably extremely difficult, unless the module is multiplicity-free or at least distributive. 
A survey of results in the   distributive case was given in 
\cite{M17}. For a Verma module $N$ over a simple Lie algebra of rank 2, $L(N)$ was determined by Jantzen \cite{J1} 3.18, using his sum formula. 
In the case that $N$ has a so-called dominant, integral highest weight, $N$ is multiplicity-free with one composition factor for every element in the Weyl group. See \cite{M} Exercise 10.5.2 for the case of $\fsl(3)$. 
To use the method suggested in the  abstract effectively, we need to find submodules $A_i$ having a simpler submodule structure than $M$. This does not help in Jantzen's examples. The examples from  Section \ref{uv2} suggest that when the method can be successfully applied, the submodule structure of $M$ turns out to be much more complex. 
 Hopefully the result mentioned in the abstract can be applied to other cases apart from those described there.
\\ \\
In her thesis 
 America Masaros, \cite{Mas} studied Verma modules $M(\lambda)$  for   the  simple Lie superalgebra 
$\fg =\osp(3,2)$ and  determined their lattice of submodules  in many cases. 
However for the Vermas $M(\gl_n)$, the lattice structure remained  elusive.  Some partial information was obtained, including the determination of all composition factors.  These modules are multiplicity free of length 8.
\\ \\
Since lattices and Lie superalgebras might be considered as distantly related 
areas, we give more detailed background on both subjects than we would otherwise. This is done in the rest of this Section and the next.  After that  we work towards finding a set of submodules 
$A_1,\ldots , A_n$ satisfying the conditions in the abstract.  
By combining results from  lattice theory and \cite{Mas} we arrive at Corollary \ref{lx} which gives  us two possibilities.  
\\ \\
To make further progress, we require a new ingredient which may be of independent interest. Namely we use 
 an important property of Sapovalov elements:   if $\gl \in \cH_\gc,$ then $ \theta^2_{\gamma} M(\gl) = 0$ (the notation is defined before Lemma \ref{k9}). This leads to a complex of Verma modules with highest weights in $\cH_\gc,$ where every map is multiplication by $\gth_\gc.$
In the present context this complex is  
$${\bf B}^\bullet: M(\gl_1)  \lra M(\gl_2)\lra M(\gl_2) \ldots$$ 
from \eqref{ec}, and a key result Corollary \ref{k6}, is that the complex is exact at $M(\gl_n)$ for $n\ge 2$.  This suggests a connection between exactness of the more general complex from \cite{M21} (5.1) and representation theory.  
Based on results from Sections \ref{elc} and \ref{osp1}, the main result for 
the $\osp(3,2)$ 
Verma modules $M(\gl_n)$  is proved in Section \ref{uv2} by building the Hasse diagram in a series of steps.
\subsection{Basic Notions of Lattice Theory.} \label{a1}
Throughout $[n]$ denotes the set of the first $n$ positive integers.  
In this paper all lattices and posets are finite. In a poset $P$, 
if $B < A$ and $B \le C \le A$ implies that $C = A$ or $C = B$ we say that $A$ {\it covers}
$B.$  The {\it 
Hasse diagram} of $P$ is a graph whose vertices are the elements of $P,$ with an edge drawn upwards from 
$B $ to $ A$ whenever $A$ {covers}
$B.$ 
A {\it lattice} is a poset $(\gL, \le)$ in which any two elements $A, B$ have a 
least upper bound (lub) $A \vee B,$ and a 
greatest lower bound (glb) $A \wedge B$ 
called the {\it  join} and {\it meet} of $A, B$ respectively.
These are necessarily unique.
Furthermore the conditions $B \le A, B = A \wedge B$ and $A = A \vee B$ are equivalent.
In this case we say that  $[B,A]$ is an {\it interval.}
 The {\it dual lattice} to $\gL$ has as 
underlying set a copy $\overline \gL$ of  $\gL$, with a bijection $X \lra \overline X$ from $\gL$ to $\overline \gL$, and then the partial order, join and meet in $\overline \gL$ are defined by $\overline X \le \overline Y,$ iff $X \ge Y,$ $\overline X \vee \overline Y = \overline{X \wedge Y}$ and 
$\overline X \wedge \overline Y = \overline{X \vee Y}$.
\\ \\
We recall some terminology from \cite{St}. 
We say that $A$ is {\it join irreducible} if $A = X \vee Y$ implies that either $A = X$ or $A = Y$. 
The glb of all elements in the lattice $\gL$, denoted by 0, is clearly join irreducible. Similarly  the lub of all elements in  $\gL$ is denoted 1. An element of that covers 0 is called an {\it atom}. Dually  an element covered by 1 is called a {\it coatom}. 
  A lattice $\gL$  
is {\it graded} if any two maximal (saturated) chains in $\gL$  have the same length $n$, and in this case we say  $\gL$ is {\it graded of degree $n$.}
\ff{In  \cite{St} such a lattice is called {\it graded of rank $n$.} But the 
extended free distributive lattice of rank $n$ 
has degree $2^n$.  So to avoid a clash of notation, we depart from the usual terminology here.
} For such a lattice, define a {\it degree function }  $\gr:\gL \lra \{0,1, \ldots n \}$ such that $\gr(0) =0$, and $\gr(A) =\gr(B) +1$ if $A$ covers $B$.  We say $\gL$ is {\it upper semimodular} if  $\gL$ is graded and 
$$\gr(B) +\gr(A) \ge \gr(B \vee A) +\gr(B \wedge A),$$ for all  $B,A\in \gL$. There is an  equivalent condition in 
\cite{St} Proposition 3.3.2. If the dual lattice  $\bar \gL$ is also upper semimodular we say  $\gL$ is {\it modular}.  (Thus if a lattice is modular, so is its dual.) In this case we have
\be \gr(B) +\gr(A) = \gr(A \vee B) +\gr(A \wedge B).\nn\ee
Alternatively 
\be \label{loy}\gr(B) -  \gr(A \wedge B)
= 
\gr(A \vee B) -\gr(A).\ee
The interval $[A, B]$ has length $\gr(B) -\gr(A)$.
According to 
\cite{B} Corollary to Theorem 15, page 41,
modularity is equivalent to the condition
$$B \vee (A \wedge
 C) = (B \vee  A)  \wedge C,$$ 
for all  $A,B, C\in \gL$ with $B\le C$. 
We say  $\gL$ is  {\it distributive} if
$$B \vee (A \wedge
 C) = (B \vee  A)  \wedge (B\vee C),$$ 
for all  $A, B, C\in \gL$.  This is equivalent to the dual condition.
$$B \wedge (A \vee 
 C) = (B \wedge A)  \vee  (B\wedge C),$$ 
for all  $A, B, C \in \gL$, \cite{DH} Theorem 6.5.3. Any 
distributive lattice is
modular, \cite{DH} Lemma 6.5.1.  
 {\it }
A sublattice of a lattice $\gL$ is a subset which is closed under the meet and join opertions.
\bexa \label{e4}{\rm Here are the  Hasse diagrams for two important lattices.  Each is minimal with respect to not satisfying certain desirable properties of lattices.    
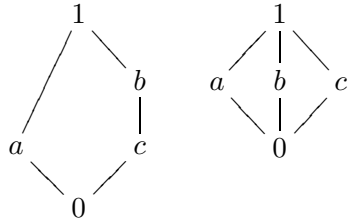
\[
\xymatrix@C=1pc@R=1pc{
& 1 &\\
&&
b \ar@{-}[ul]&\\
a \ar@{-}[dr] \ar@{-}[uur] &&
c \ar@{-}[dl] \ar@{-}[u]&\\
& 0 &}
\xymatrix@C=1pc@R=1pc{
& 1 &\\
a \ar@{-}[dr] \ar@{-}[ur] &b \ar@{-}[d] \ar@{-}[u]&
c \ar@{-}[dl] \ar@{-}[ul]&\\
& 0 &}
\] Indeed a lattice is modular (resp. distributive) iff it contains no sublattice with Hasse diagram as in the first (resp. both) diagrams above, 
\cite{LP} Theorem 1.24 
(resp. Theorem 1.29). 
The second diagram above is the lattice of subgroups of the Klein four group $V$.  As a $\Z$-module $V$ is semisimple of length two, and has a unique composition factor up to isomorphism.}\eexa \noi
\subsection{Modular Lattices.} \label{u1}
We sometimes refer to elements of a lattice 
 $\gL$  as submodules.  
First we recall some definitions from  \cite{M17}. Let $\equiv$ be the smallest equivalence relation on the set of intervals such that
\be \label{e2} 
[ A \wedge B,B]
\equiv [A,A \vee B]\ee 
We sometimes write the interval $[X, Y ]$ as $Y/X$. Thus \eqref{e2} is motivated by the second isomorphism theorem:
 $B/(A\cap B)\cong (A + B)/A.$ The  equivalence relation $\equiv$ is not well behaved for non-modular lattices.  Indeed from the first diagram in Example \ref{e4} we have since $1=a \vee b$ and $0 = a  \wedge  b$, that $1/a   \equiv b/0.$ 
Thus an interval of length 1 can be equivalent to an interval of length  2. However for modular lattices, equivalent intervals have the same length by 
\eqref{loy}. For a module $M$ the lattice of submodules  $L(M)$ of $M$  is modular by Dedekind's modular law, see for example \cite{M} Lemma 1.2.6.  From now on we assume that all lattices are modular. We call an interval of length one $[X, Y ]$ or $Y/X$, a {\it simple lattice factor} of  $\gL$.
The equivalence classes under $\equiv$ restricted to the set of
intervals of length one 
will be called {\it lattice composition factors} of  $\gL$. If $J\neq 0$ is join irreducible it has a unique maximal submodule, denoted by
 $J^0$. Define $L_J=J/J^0.$
\bl \label{L8} 
Let $\gL$ be a modular lattice. If $M/N$  is a {simple lattice factor} of  
$\gL$, then  $M/N\equiv J/J^0$ where $J\neq 0$ is join irreducible in  
$\gL$ and $M \ge J$. \el 
\bpf If the result is false choose a counterexample with $M$ minimal under the partial order $\ge$.  Then $M$ is not join irreducible, so $N$ is not the only maximal submodule of $M$. 
 Let $P$ be another  maximal submodule. Then $M = N \vee P$, so
$M/N\equiv P/(P \wedge N)$ with $M > P$. By the minimality of $M$,  $P/(P \wedge N) \equiv J/J^0$ with  $J$  join irreducible. This contradiction yields the result.
\epf  
\noi 
Consider  a maximal chain in the lattice $\gL$:
\be \label{sow}0=X_0 < X_1 < \ldots < X_n =1.\ee We say the lattice composition factor $a$ has multiplicity $k$ in 
\eqref{sow} if exactly  $k$ of the intervals $[X_{i-1},X_i]$, $1 \le i \le n$ belong to the equivalence class $a$.
There is  Jordan-H$\ddot{o}$lder Theorem for modular lattices. 
  \bt\label{bat} In a modular lattice $\gL$, the multiplicities of a lattice composition factor in any two maximal chains are equal. \et
\bpf If the result is false, let $\gL$ be a counterexample of minimal degree $n$. Since $\gL$  is modular, all maximal chains in
$\gL$ have length $n$. If 
$\gL$ had a unique atom $A$, then  the result holds for the interval $[A,1]$, and any maximal chain would start with $0<A$.  Thus $\gL$ would not be a counterexample.  So suppose $A, B$ are distinct atoms. Then $A\wedge B=0$, so  by modularity $A\vee B$ covers both $A$ and $B$.
  Let $A\vee B=Y_2 < \ldots < Y_n =1$ be a maximal chain in the interval  $[A\vee B, 1]$.  Since the interval $[A,1]$ is  not a counterexample, any maximal chain in $[A,1]$ 
 has the same multiplicities as the chain $$A<A\vee B=Y_2 < \ldots < Y_n =1.$$ 
Thus any maximal chain in $\gL$ 
that begins with $0 < A$ has the same multiplicities as the chain
 \be \label{lay}  0<A<A\vee B=Y_2 < \ldots < Y_n =1.\ee
But then reversing the roles of $A$ and $B$ and using the definition of equivalence,  shows
that any maximal chain that begins with $0 <B$ has the same multiplicities as the chain in \eqref{lay}. Since any maximal chain must start with an atom, and  $A, B$  are arbitrary, $\gL$ is  not in fact a counterexample.
\epf\noi 
We
say that $\gL$ is {\it multiplicity-free} if 
the multiplicity of any  lattice composition factor in a maximal chain
equals 1. Similarly a modue $M$ is {\it multiplicity-free} if 
the multiplicity of  any composition factor equals 1.
\subsection{Distributive Lattices.} \label{cat}
 For a module $M$  the multiplicity-free condition on $M$  is closely related to distributivity of $L(M)$.  It is well known that $L(M)$ is distributive iff every semisimple subfactor of  $M$   is multiplicity-free. It is hard to track down the first citation, but a short proof is given in \cite{M17}  Proposition 2.1. According to  Stephenson 
\cite{S} Proposition 1.1, $L(M)$ is distributive iff for every  subfactor of $M$ which is a direct sum $A \op B$ of two non-zero submodules, we  have $\Hom(A,B) =0$. It is not hard to deduce the result from this. 
\\ \\
We say a lattice is {\it restricted} if it has at least 2 atoms and 2 coatoms.  A restricted  lattice 
$\gL$ has {\it rank} if $n$ is minimal  such that $\gL$ can be generated by $n$ join irreducibles $J_1,\ldots , J_n$ using the meet and join operations.    For 
$\emptyset \neq I \subseteq [n] = \{ 1, \ldots, n\}$, set  \ff{If $X, Y$ are disjoint subsets (of some multiplicative abelian group), and $\prod_{X}$ denotes the product of elements in $X$, then $\prod_{X\cup Y} = \prod_{X} \prod_{Y}$.  Taking $X$ to be empty tells us that the empty product equals 1, and similarly  the empty 
 sum is zero.  But the same argument applied to $\bigvee_{I} J$ and  $\bigwedge_{I} J$ tells us only that $\bigvee_{\emptyset} J \le \bigvee_{i\in I} J_i \mbox{ and  }\bigwedge_{I} J \le  \bigwedge_{\emptyset} J$ for all $I \subseteq [n]$.  For extended $\gL_n$ (defined in Subsection \ref{dog}) this gives two possibilities for 
$\bigvee_{\emptyset} J   \mbox{ and  }\bigwedge_{\emptyset} J$. For this reason we prefer not to think about empty meets and joins.
}
\be \label{t1} \bigvee_{I} J = \bigvee_{i\in I} J_i \mbox{ and  }\bigwedge_{I} J = \bigwedge_{i\in I} J_i.\ee

\bl \label{asp} 
Each $\bigwedge_{I} J$ is join irreducible. 
\el
\bpf  
By induction on $n $. For $n>1,$ the previous step  
takes care of all cases where $|I|\le n-1$. In the remaining case
 $\bigwedge_{[n]} J$ is the unique minimal element of the lattice.
\epf 
\bexa {\rm The power set $\B_n$ on $[n]$ is a distributive lattice with union as join, and intersection as meet. We call $\B_n$ the Boolean lattice (or poset) of rank $n$. We also call $\B_3$  {\it the box lattice} $\B$ because of its Hasse diagram, \cite{La} Figure 5.2. The lattice $\B_n$ 
is generated by the join irreducibles $\{1\}, \{2\}, \ldots \{n\}$ and has $n$ composition factors up to isomorphism. For  example $\B$ has  3 composition factors up to isomorphism, because parallel edges of the box are equivalent.   
}\eexa
\noi 
If $P$ is a poset, a subset $I$ of $P$ is called a {\it down-set} (or {\it order ideal}) if $x \in I, y \in P$ and $y\le x$ implies $y \in I.$ 
If $x \in P$, the down-set $\gL_x= \{y\in P| y \le x\}$ 
is known as a {\it principal down-set}. 
Denote the set of down-sets of $P$ by $J(P)$. Then $J(P)$ is a distributive lattice taking union and intersection as $ \vee, \wedge$ respectively. Every finite distributive lattice  has this form.  This is the content of the next result known as the fundamental theorem on finite distributive lattices.
 
\bt \label{t4} Let $L$ be a finite distributive lattice. Then there is a unique, up to isomorphism, finite poset $P$ such that  $L \cong J(P)$\et
\bpf  Given a lattice $L$ as above, let $P$ be the set of join  irreducible elements of $L$. This is a poset with the order induced from that on $L$. It is shown in \cite{St} Theorem 3.4.1 that $L\cong J(P)$ as  lattices. See also \cite{B} Theorem 3, page 59.
 \epf \noi 
When  $P=\B_n$, $J(P)$ is isomorphic to extended $\gL_n$ which has $2^n$  composition factors up to isomorphism. This could be an argument in favor of defining the free distributive lattice 
 of rank $n$ to be extended $\gL_n$.
\\ \\
A method for drawing the Hasse diagram of   $J(P)$ given $P,$ is given in 
\cite{St} pages 108-109, follwed by an example where $P$ is a zig-zag poset.  We use this method in Section   \ref{uv2}, to determine the lattice of submodules of our  Verma modules. We briefly recall the details of Stanley's method.  The idea is that $J(P)$ can be constructed using $P$ by gluing together Boolean lattices.  Begin with 
the set $I$ of minimal elements of $P$. If $|I|=m$ form the lattice $J(I)\cong \B_m.$ If $I\neq P$ choose a minimal element $x$ of $P  \backslash I$. Adjoin a join irreducible to 
$J(I)$ covering the down set
$\gL_x  \backslash \{x\}$. The set of joins of elements covering 
$\gL_x  \backslash \{x\}$ forms a Boolean algebra, draw in any new joins required to show this. (An example occurs passing from Fig. 3-13 to Fig. 3-14 in \cite{St}.) There still may be elements that do not yet have joins. Add in these joins and repeat until we obtain $J(P)$. (An example occurs passing from Fig. 3-19 to Fig. 3-21 in \cite{St}.)

\subsection{The Free Distributive Lattice of Rank $n$.} \label{dog}
A {\it proposition} (or Boolean variable) is a  variable that can take on the values T or F. Propositions $P_1, \ldots, P_n$ are {\it independent} if the values of  the $P_i$ can be assigned independently of each other.  
For propositions $X, Y$, define $X \wedge 
 Y$, (resp. $X \vee Y)$ by requiring that  $X \wedge Y = $ T iff $X =$ T and $ Y= $ T
 (resp. $X \vee Y = $ T iff $X =$ T or $ Y= $ T). Then define $X \le Y$ to mean $X \wedge Y= X$. Thus $X \le Y$ means that $X$ implies $Y$. 
The  {\it free distributive lattice} $\gL_n$ of rank $n$ is the lattice of propositions generated from $P_1, \ldots, P_n$ by using the join  $\vee$ and  meet $\wedge$ operations.  
We define  $ \bigvee_{I} P$ and $\bigwedge_{I} P,$ by analogy with \eqref{t1}.
There is no universal agreement on the definition of $\gL_n$.  Our definition is equivalent to that given in \cite{DH} Section 6.8, where a labelled diagram of the lattice $\gL_3$ can be found, see also Section \ref{uv2} and \cite{B} Figure 8, page 33.   We obtain {\it extended } $\gL_n$ from $\gL_n$ by adjoining additional elements $\hat 0$ and   $\hat 1$, such that  $\hat 0 <X$ and $X< \hat 1$ for all $X\in \gL_n$. Often the 
free distributive lattice of rank $n$ is defined as  extended  $\gL_n$. For 
additional clarity we sometimes refer to $\gL_n$ as {\it restricted } $\gL_n$.
\br \label{L1}{\rm The lattice $\gL_n$ is self dual. Indeed, if  $P_1, \ldots, P_n$ are independent propositions,  then so too are their negations. So the statement follows from DeMorgan's Laws.} \er \noi 
The Wikipedia article on Dedekind numbers \ff{
$https://en.wikipedia.org/wiki/Dedekind\_number$} gives the 
  Hasse diagram of extended $\gL_3$.  The diagram is labelled using independent Boolean variables $A, B, C$. The number of elements of extended $\gL_n$ is called the $n^{th}$ {\it Dedekind number} $M(n)$, and the values of   $M(n)$ are only known for $n\le 8$.  For example $M(3) = 20$ and $M(8)$  is approximately 
 $5.61 \times 10^{22}$.  
We could define a function $f$ defined on $\N$ to have exponential growth if asymptotically $\log_2 f(n)$ grows  like a linear polynomial in $n$.  
According to Wikipedia,
$$\left( \begin{array}{c} n\\ \lfloor n/2\rfloor\end{array}\right) 
\le \log_2 M(n)\le \left( \begin{array}{c}n\\ \lfloor n/2\rfloor\end{array}\right)(1 +O\left( \begin{array}{c}\log  n \\ \hline n \end{array}\right)
) .$$ 
Now $\left( \begin{array}{c} n\\ \lfloor n/2\rfloor\end{array}\right) $ is a polynomial in $n$ of degree $ \lfloor n/2\rfloor$.  As $n$ increases, this means that $\log_2 M(n)$ eventually grows faster than any polynomial in $n$.  
Thus the growth of $M(n)$ is much faster than exponential. \\ \\
  In the  Wikipedia article several equivalent definitions are given for extended $\gL_n$, one of which involves involves Boolean variables.  In the diagram for $\gL_3$ , 0 an 1 are labelled as ``contradiction" and ``tautology" respectively.  However we don't regard these as genuine propositions, since they can only take on a single value.  Also  $P_1, \ldots, P_n$  already have a lub and glb of $\bigvee_{[n]} P $ and  $\bigwedge_{[n]} P$  respectively. \\ \\
We briefly explain the 
theory of Disjunctive Normal Form (DNF), following  \cite{La}.  Begin with independent Boolean variables $P_1, \ldots, P_n$ as above. A {\it literal} is either one of the $P_i$ or its negation. A {\it Boolean function} is a proposition that can be 
formed from the literals using meets and joins. A
{\it miniterm} is a conjuction, that is a meet of literals.  Any Boolean function has an essentially unique expression (called the DNF) as an irredundant  disjunction or join of miniterms \cite{La} Proposition 5.9 and Theorem 5.10.  
More precisely, the uniqueness is up to changing the order of the  literals in miniterms, and changing the order of the   miniterms in the disjunction. DNF is useful in efficient circuit design. If we use {\it and gates} and {\it or  gates} that admit multiple inputs, then we need only a single {\it or  gate}, and an {\it and gate} for each miniterm, as well as a number or {\it not  gates}.  DNF applies to the lattice    
$\gL_n$  provided we redefine a literal to be one of the $P_i$. Negation does not preserve $\gL_n$.  Thus for a  Boolean function in $\gL_n$, no {\it not  gates} are needed. \noi 
\bl \label{L2} \bi
\itema $\bigwedge_{I} P = \bigwedge_{J} P$ iff $I=J$. 
\itemb 
Any join irreducible element of 
$\gL_n$ is equal to $\bigwedge_{I} P$ for some $I$.
 \ei\el
\bpf Both parts  follow from the uniqueness of disjunctive normal form. 
An easier proof of (a) is as follows.  Suppose $i\in I,$ but $i\notin J$. Assign values $P_i: =$ F, and $P_j: =$ T for all $j\in J.$  Then  $\bigwedge_{I} P =$ F and $ \bigwedge_{J} P$ =T.
\epf 

\br {\rm The {\it free modular lattice} $M_n$ of rank $n$ was first considered by Dedekind \cite{D}, who showed that  $M_3$ has 28 elements.  A proof can also be found in \cite{B}   Section III.6, where it is also shown that $M_4$ is infinite.  In \cite{F} it is shown that $M_5$ has unsolvable word problem.}  
\er
\bt \label{t2} Suppose $\gL$ is  a restricted distributive lattice of rank $n$ generated by  join irreducibles $J_1,\ldots , J_n$.  Then there is a surjective map of lattices restricted 
$\gL_n \lra \gL$ sending $P_i$ to $J_i$.\et 
\bpf See \cite{B} Ch III.4.\epf \noi 
\br {\rm 
The analog of 
Lemma \ref{L2} (a) fails for the box lattice $\B$ because $\{1\} \cap \{2\}= \{1\}\cap\{3\}= \{2\}\cap\{3\}.$ It is an interesting exercise to construct a surjective map of lattices from $\gL_3$ to $\B$, and check that it is well defined. 
}\er

\bexa \label{e8}{\rm Below  is the lattice of subgroups of the cyclic group $\Z/(12)$ of order 12.  Note that $\Z/(12)$ is not multiplicity free but $L(\Z/(12))$ is distributive.  The smallest module with this property is the cyclic group of order 4, but $L(\Z/(12))$ will reappear in Lemma \ref{hag}.
\[
\xymatrix@C=1pc@R=1pc{
&\Z/(12)&\\
2\Z/(12)\ar@{-}[dr] \ar@{-}[ur] &&
3\Z/(12) \ar@{-}[dl]  
\ar@{-}[dr]\ar@{-}[ul]&\\
& 4\Z/(12)& &6\Z/(12)\ar@{-}[dl] \\
 &&
0 \ar@{-}[ul]&\\
} 
\] }
\eexa
\noi The module $\Z/(12)$ contains isomorphic simple subfactors which do not isomorphic correspond to simple lattice factors.   Their isomorphism is not a consequence of the second isomorphism  theorem.  We call such isomorphisms {\it acccidental}.
\bl \label{cow} If the module $M$ is muliplicity free then $L(M)$ has no accidental 
isomorphisms. In other words, if two simple subquotients of $M$ are isomorphic as modules, then the corresponding intervals are equivalent in $L(M)$. \el
\bpf If the result is false, then by Lemma \ref{L8}, $M$ contains distinct join irreducible submodules $J,K$ such that $J/J^0 \cong K/K^0$.  
If $K  \subseteq J,$ then $M$ contains the chain of submodules
$$K^0 \subset  K  \subseteq J^0 \subset J,$$
 with two isomorphic composition factors, a contradiction.  Similarly $K$ cannot contain $J.$ Thus $J \cap K $ is a proper submodule of both $J$ and $K.$
Since $J^0, K^0$ are the unique maximal submodules of $J, K$ it follows that $J \cap K  \subseteq J^0 \cap K^0$. Thus $J + K /J \cap K $ is a subfactor of $M$ that is also a counterexample. Starting again with  $M= J + K $ and $J \cap K =0$ we see that $M= J\op K$.  But then $M$ is not  muliplicity free. 
\epf
\bc \label{hog} If the module $M$ is muliplicity free, then the map $J \lra L_J$ induces a bijection between join irreducibles $J ($with $J\neq 0)$ in $L(M)$ and isomorphism classes of composition factors of $M$.
\ec
\bpf Combine Lemma \ref{L8} and Lemma \ref{cow}. \epf
\noi We record some well known examples of lattices from \cite{St} Example 3.1.1. The lattice  $[{\bf n}]$ is the set ${[ n]}$ with its usual order.  This lattice is {\it uniserial}, that is any two elements are comparable.  If $n$ is a positive integer, then the lattice of subgroups of  $\Z/(n)$ is denoted $D_{ n}$. The following result is well known.
\bl  \label{hag} The list of distributive lattices of length 3 is given as follows.  If the lattice does not have a name, we give its Hasse diagram.
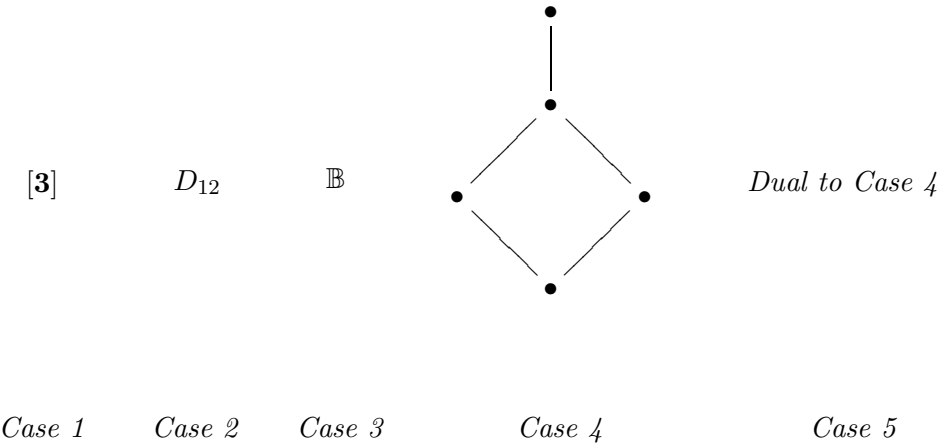
\[ 
\xymatrix{
\\ \\
\mbox{
{\rm[{\bf  3}]}
}
\\ \\
}
\quad \quad\quad
\xymatrix{
\\ \\
D_{12}
\\ \\
}\quad\quad\quad
\xymatrix{
\\ \\
\B
\\ \\
}\quad\quad\quad
\xymatrix{
&
\bullet\ar@{-}[d]&\\
&
\bullet&\\
\bullet \ar@{-}[dr] \ar@{-}[ur]
 &&
\bullet \ar@{-}[dl] \ar@{-}[ul]&\\
& 
\bullet  & \\
&
&}
\xymatrix{
\\ \\
\mbox{Dual to Case 4}
\\ \\}\] 
\[ 
\mbox{Case 1}
\quad\quad\;
\mbox{Case 2}\quad\quad
\mbox{Case 3}\quad\quad\quad\quad\;\;
\mbox{Case 4}\quad\quad\quad\quad\quad\quad\quad
\mbox{Case 5}\quad\quad\] 
\el
\bpf We work throughout up  to lattice isomorphism. Let $\gL$ be a lattice as in the statement of the Lemma,  and denote the lub of its elements by 1.  If 1 is join irreducible, it is easy to see we have Case 1 or 4.  Now suppose $\gL$  has two (distinct) coatoms $b, d$.  Then $1=  b\vee d$.  Clearly $b\wedge d$ covers 0. If both    $b, d$ are join irreducible,  we have Case 5.  
If exactly one is join  irreducible, we have Case 2. If neither is join  irreducible,  we obtain the lattice with Hasse diagram
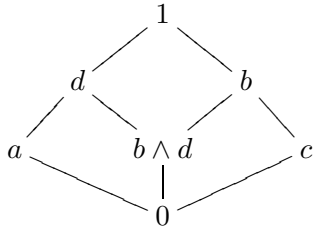
\[\xymatrix@C=1pc@R=1pc{
&& 1 &\\
&d \ar@{-}[dr]\ar@{-}[dl]
 \ar@{-}[ur] &
&
b
\ar@{-}[ul]
&\\a \ar@{-}[drr]&&
b\wedge d\ar@{-}[d]
\ar@{-}[ur] &&
c\ar@{-}[dll] \ar@{-}[ul]&\\
&& 0 &&}
\] 
But then the sublattice consisting of the elements $0, 1, a, b, c$ is the 5 element non-modular  lattice from Example \ref{e4}. So $\gL$ is not distributive.  
If  $\gL$  has $n \ge 3$ coatoms, then the dual lattice would contain $n$ atoms, and hence contain the self-dual sublattice $\B_n$ which has length $n$. Thus the only possibility is $n=3$ and we have Case 3.
\epf
\brs{ \rm The lattices listed in the Lemma have the form $J(P)$ for the five posets $P$ of size three listed in \cite{St} Figure 3-1, but it will be useful to know their Hasse diagrams. Clearly there is a bijection between the lattices $\gL$ in the Lemma and  distributive lattices of length 4 with a unique atom. If $\gL$ is as in the Lemma, we call the corresponding length 4 lattice {\it augmented} $\gL$.   }\ers
\br{\rm We outline how to prove the statement in the abstract, that we can reconstruct $L(M)$ from a set of submodules $A_1,\ldots , A_n$  which contain all composition factors of $M$, and such that each $L(A_i)$ is known.  
Using Lemma \ref{L8}, we can assume the $A_i$ 
are join irreducble.  This may entail changing the original set of submodules, see Subsection \ref{lq}.  From this information, by Corollary \ref{hog}, we know the structure of the poset of join irreducibles, and then we can use Stanley's method.}\er
\section{Lie Superalgebras.} \label{osp}

\subsection{Introduction.} \label{b1}
We use \cite{H}, \cite{K1}  and  \cite{K}, \cite{M}  as general  references for
Lie algebras and 
Lie superalgebras respectively. A {\it Lie superalgebra} is a $\Z_2$-graded vector space 
$\fg=\fg_0 \op \fg_1$  together with a bilinear map $[\;,\;]:\fg \ti \fg \lra \fg$ satisfying a $\Z_2$-graded version of the Jacobi identity and skew-symmetry, \cite{M}  Section 1.1.  The $\Z_2$-grading on 
$\fg$ means that $[\fg_i, \fg_j]\subseteq \fg_{i+j}$ for $i,j \in \Z_2$. This means that  $\fg_0$ is a Lie algebra and $\fg_1$  is a $\fg_0$-module via 
$[\;,\;] $. Finite dimensional simple 
Lie superalgebras over an algebraically closed  field $\ttk$ of characteristic zero were classified by Kac. Several special cases of the classification were obtained by other authors, see \cite{K} pages 47-48 for details.  Such an algebra $\fg$ is called {\it classical} if $\fg_0$ is reductive.  
\ff{A Lie algebra $\fk$ is {\it reductive} if $[\fk,\fk]$ is semisimple, that is a direct sum of simple Lie algebras.  If $\fk$ is { reductive} then $\fk=\fz(\fk) \op [\fk,\fk]$, where 
$\fz(\fk)= \{x \in \fk| [x,y] = 0 \mbox{ for all } y\in \fk\}$ is the {\it center}  of $\fk.$}
The remainder are called {\it of Cartan type}.  From now on we consider only    finite dimensional classical simple Lie superalgebras $\fg$  over $\ttk$.  The classification of such 
Lie superalgebras can be found in \cite{M} Theorem 1.3.1.
\\ \\
For simplicity we assume that $\fg \neq \fp(n), \fq(n)$, two  infinite families  of algebras in the Theorem just cited. The remaining Lie superalgebras are called {\it basic classical simple}. Let $\fh$ be a Cartan subalgebra of $\fg_0$.  This is an {\it abelian} subalgebra (meaning that $[\fh,\fh]=0$) of $\fg_0$ such that \eqref{b9} below holds. 
For $\ga\in\fh^*$, set 
\be \label{b7}  \fg^\ga =\{x \in \fg| [h,x] = \ga(h)x\}.\ee We say that $\ga$ is a 
{\it root of $\fg$} if $\ga \neq 0$ and $\fg^\ga \neq 0$. Denote by the set of roots of $\fg$ by $\Delta.$ The {\it adjoint action} $(\fh,\fg)\lra \fg$ of $\fh$ on $\fg$, defined by 
$ (h,x)\lra [h,x]$ is diagonalizable and clearly each $\fg^\ga$ is an eigenspace for this action. We have a   {\it root space decomposition}, \cite{H}  Chapter 8, \cite{M} 2.1, 8.1,
\be \label{b9} \fg=\fh\op\bigoplus_{\ga\in\Gd}\fg^\ga.
\ee
Moreover if $\fg\neq\fpsl(2,2)$ (which we will assume for simplicity), then each root space in \eqref{b9} has dimension 1 over $\ttk$ and we choose a basis element $e_\ga$ for $\fg^\ga$. 
A root $\ga$ is called {\it even} (resp. {\it odd}) if $\fg^\ga \subset \fg_0$ (resp. $\fg^\ga \subset \fg_1$). Denote the set of even (resp. odd) roots by $\Delta_0$
(resp.  $\Delta_1$). 
We have a disjoint union $\Delta=\Delta_0 \cup\Delta_1.$ This fails for $\fq(n)$.
There is a similar  { root space decomposition} for $\fg_0$,  
\be \label{b91} \fg_0=\fh\op\bigoplus_{\ga\in\Gd_0}\fg^\ga.\ee
In addition there is  a disjoint union $\Delta=\Delta^+ \cup -\Delta^+,$ with the property that $\ga, \gb\in \Delta^+, \ga + \gb\in \Delta$ implies $\ga + \gb\in \Delta^+$. This fails for $\fp(n)$.
For $i=0,1$, set  $\Delta_i^{\pm}= \Delta_i\cap \Delta^{\pm}.$ 
Define $\fn^\pm =\bigoplus_{\ga\in\Gd^\pm}\fg^\ga.$ From the defining properties of $\Gd^\pm$ and \eqref{b9}, it follows that $\fn^\pm$ are subalgebras of $\fg$  and  we have a {\it a triangular decomposition} of $\fg.$
\be \label{gtri}\mathfrak{g} = \mathfrak{n}^- \oplus \mathfrak{h}
\oplus \mathfrak{n}^+.\ee
There is  a similar  {triangular decomposition} of $\fg_0$ with $\fn_0^\pm =\bigoplus_{\ga\in\Gd_0^\pm}\fg^\ga,$
\be \label{stri}\mathfrak{g}_0 = \mathfrak{n}_0^- \oplus \mathfrak{h}
\oplus \mathfrak{n}_0^+.\ee
A root in $\Delta^+ $ is {\it simple} if it cannot be written as a sum of two other roots in $\Delta^+ $.  Denote the set of simple roots  by $\Pi$. In  \cite{H}  Chapter 10,  $\Delta^+ $ and $\Pi$ are denoted by  $\Phi^+ $ and  $\Delta $ respectively. For the Lie superalgebra versions, see  \cite{M} 3.4 and 8.1.
We use the {\it Borel subalgebras} $\mathfrak{b} =  \mathfrak{h}
\oplus \mathfrak{n}^+$ and
$\fb^- =
\mathfrak{n}^- \oplus \mathfrak{h}$.
Equation \eqref{gtri} leads to some other useful decompositions
\be \label{gto} \mathfrak{g} = \mathfrak{n}^- 
\oplus \mathfrak{b}, \quad \mathfrak{b} = 
\mathfrak{h}
\oplus \mathfrak{n}^+.\ee
Since $\fn^\pm$ is defined as a sum of root spaces, it follows that 
 $[\fh, \fn^\pm] \subseteq \fn^\pm$ (in fact equality holds). Then from the  second equation in \eqref{gto}, we have 
\be \label{gtp}  [\fb, \fn^+] \subseteq \fn^+ .\ee  This says that
$\fn^+$ is an {\it ideal} in $\fb$.   Thus there is a natural Lie superalgebra structure on  $\fb/\fn^+$, and we have  $\fb/\fn^+\cong \fh$ as algebras. 
Define 
\be \label{b3} \gr_0=1/2\sum_{\ga \in\Delta^{+}_0}\ga, \quad  \gr_1=1/2\sum_{\ga \in\Delta^{+}_1}\ga, \quad  \gr =
\gr_0-\gr_1.
\ee
\noi
There is a uniform construction (due to Kac) of basic classical simple Lie superalgebras
using contragredient Lie superalgebras, often called Kac-Moody (KM) Lie superalgebras, beginning with a generalized Cartan matrix $A$, \cite{M} Chapter 5. The KM Lie algebras of finite type are exactly the finite dimensional simple Lie algebras, 
\cite{K1},  and the simple KM Lie superalgebras of finite type are exactly the finite dimensional basic classical simple Lie superalgebras. Replacing $A$ by a $DA$ for a non-singular diagonal matrix $D$ results in an isomorphic algebra. In the finite type case this can always be done with $DA$ symmetric, and then the symmetrized Cartan matrix can be used to define an even non-degenerate invariant form on the simple algebras, \cite{K1} Theorem 2.2, \cite{M} Section 5.4 whose restriction $<\;,\;>$ to the Cartan subalgebra $\fh$ is also non-degenerate. If $A$ is non-singular, which in finite type is always true in the Lie algebra case, as it is for all basic classical simple superalgebras different from $\fpsl(n,n)$ for $n\ge 2$, then the simple roots form a basis for $\fh$. For simplicity we assume   $\fg\neq\fpsl(n,n)$. We can use the non-degeneracy of the bilinear form to define a linear  isomorphism $\fh\lra \fh^*$. The inverse isomorphism $\ga\mapsto h_\ga$ is defined by $<h_\ga,h> = \ga(h)$ for all $h \in \fh$. 
Then we can define a bilinear form $(\;,\;)$ on $\fh^*$ by 
\be \label{b6} (\ga, \gb) =<h_\ga, h_\gb>\quad  (= \ga( h_\gb) = h_\gb(\ga)).
\ee
We say that a root $\ga$ is {\it isotropic} if $(\ga,\ga)=0.$ Such roots do not exist for simple Lie algebras and any isotropic root is odd. If $\ga$ is a non-isotropic,  set $\ga^\vee = 2\ga/(\ga,\ga).$ The {\it reflection} $s_\ga: \fh^*\lra \fh^*$ is defined by $s_\ga(\gl)= \gl - (\gl,\ga^\vee)\ga$. Note that $s_\ga$ sends $\ga$ to $-\ga$ and fixes the hyperplane, $\cH = \{\gl|(\gl,\ga)=0\}.$ Thus geometrically $s_\ga$ is the  { reflection} in the hyperplane orthogonal to $\ga$. The subgroup, $W$ of End$(\fh^*)$ generated by all  reflections  is called the {\it  Weyl group}.
The {\it dot action} of the Weyl group on $\fh^*$ is defined by $(w, \gl) \lra w\cdot \gl = w(\gl+\gr) -\gr$.  This can be thought of as a version of the usual $W$-action, where the origin is shifted to $-\gr$.
All reflecting hyperplanes for the { dot action} pass through $-\gr$.

\subsection{Enveloping Algebras.} \label{b5}
The {\it enveloping algebra} of a Lie superalgebra $\fk$ is a $\Z_2$-graded associative $\ttk$-algebra $U(\fk)$  that has the same representation theory as $\fk$, (a $U(\fk)$-module is the 
same thing  as a $\fk$-module) see \cite{M} Chapter 6. This allows methods from ring theory to be use to study modules for $\fk$. 
The Poincare-Birkhoff-Witt theorem (PBW) gives a vector space basis for 
$U(\fk)$.  From PBW it follows that $\fk$ embeds as a subspace of $U(\fk)$, and the image of $\fk$ generates $U(\fk)$ as an algebra. 
As another  consequence of PBW and \eqref{gto}, we have as vector spaces 
\be \label{gta}U(\fb) = U(\fh)\ot U(\fn^+) , \mbox{ and } U (\fg) = U(\fn^-) \ot U(\fb),\ee   \cite{M} Lemma 6.1.4. Equation \eqref{gta} has some important consequences. First there is an $\ttk$-algebra map $\gve: U(\fn^+) \lra \ttk$ sending each $x\in  \fn^+$ to zero, and 
$\ker( \gve)= U(\fn^+) \fn^+= 
\fn^+U(\fn^+)$. 
\ff{$\gve$ is the counit in the bialgebra structure of $U(\fn^+)$ 
 and 
$\ker( \gve)$ is called the augmentation ideal of $U(\fn^+)$.}
From \eqref{gtp} we see that $U(\fb) \fn^+= 
\fn^+U(\fb)$  is a two sided ideal of $U(\fb)$. 
Tensoring the exact sequence 
$$ 0 \lra
U(\fn^+) \fn^+ \lra U(\fn^+){\stackrel{\gve}{{\lra}}}\ttk \lra 0$$
with $U(\fh)$ and using \eqref{gta}, we see that 
$U(\fb)/U(\fb) \fn^+\cong U(\fh)$ as $\ttk$-algebras. Thus any $U(\fh)$-module can be regarded as a $U(\fb)$-module which is annihilated by 
$U(\fb) \fn^+.$ 
 The {\it Verma module} 
$M(\gl)$ with highest weight $\gl \in \fh^*$, and highest weight vector $v_\lambda$ is defined as follows. Start with a one dimensional $U(\fh)$-module  $\ttk v_\gl$ with weight $\gl$.  By the above remarks, $\ttk v_\gl$ can be regarded as a $U(\fb)$-module with $\fn^+ v_\gl=0$. Then define
$M(\gl)=  U (\fg)\ot_{ U(\fb)} \ttk v_\gl$ as an induced module. There is a unique simple factor of $M(\gl)$, which we denote by $L(\gl)$ \cite{M} 8.2.
From \eqref{gta}, we have
$$M(\gl)=U(\fn^-) \ot U(\fb)\ot_{ U(\fb)} \ot \ttk v_\gl
= U(\fn^-) \ot \ttk v_\gl,$$
so $M(\gl)$ is a free $U(\fn^-)$-module of rank 1, generated by $v_\gl$.  
 In  \cite{H}  Chapters 20 and  21,  $M(\gl)$ and $L(\gl)$ are denoted by  $Z(\gl)$ (called a standard cyclic module) and  $V(\gl) $ respectively.

\subsection{Uniqueness up to Conjugacy of some Constructions.} \label{u2}
Let $\fg$ be a classical simple Lie superalgebra with $\fg \neq \fp(n), \fq(n)$
and set $\fk= \fg_0.$ From the definitions we see that any one of
 $\fn^+, \Delta^+ $ or  $\Pi$ determines the other two. 
A similar statement holds for $\fk$. 
By \cite{H} Theorem 16.4 and Corollary 16.4 any two Borel or Cartan subalgebras of $\fk$ are conjugate under the group  $\cE(\fk)$  of inner automorphisms. 
Also by \cite{H} Theorem 10.3, any two sets of simple roots of $\fk$ are conjugate under $W$.  Hence the definitions of $\fh,\fb,\fn^+, \Delta^+ $, $\Pi, \gr$, the dot action  
and Verma modules could be regarded as canonical in the reductive case.
\\ \\
For $\fg$ the situation is more complicated, though well understood.  Any element of $\cE(\fk)$ 
extends uniquely to an  automorphism
of  $\fg$. Thus concerning Borel subalgebras, we can fix $\fh$ and a Borel {\bf b} in 
$\fk$, and then any Borel in $\fg$ 
is conjugate to another Borel $\fb$ with 
$\fb_0=  {\bf b}$.  Now by the proof of 
\cite{M} Theorem 3.1.2, there are only finitely many Borels $\fb$ of 
$\fg$ with $\fb_0=  {\bf b}$. 
For $\osp(3,2)$ there are exactly two such.  This means that in the next Subsection, we will need two versions of   $\Pi, \Gd^+, \fb,\gr, M(\gl)$ and the dot action.  The relation between the two versions of $\Delta^+ $ is easily explained: we can obtain one from the other by changing the sign of an isotropic simple root. Apart from this the two versions have the same positive roots.  This operation or its extension to $\fb$ is sometimes called an {\it odd reflection}.  More generally, by  \cite{M} Theorem 3.1.3, any two Borels in 
$\fg$ with the same even part are related by  a sequence of  { odd reflections}.
\subsection{The Lie Superalgebra $\osp(3,2)$.} \label{b2}
Now suppose $\fg$ is the classical simple Lie superalgebra $\osp(3,2)$. 
In this case $\fg_0\cong \fsl(2) \ti \fsl(2) $, and  as positive roots for  
$\fg_0$ we take $\Delta^{+}_0=\{\gep, 2\gd\}$. Then there are 2 choices of roots $\Delta^{+}_1$ of $\fg_1$, such that $\Delta^{+}=\Delta^{+}_0 \cup \Delta^{+}_1$ is a set of positive roots for  
$\fg$. These sets and the corresponding set of simple roots $\Pi$ for $\fg$ are given by
\[ \Pi^{\rm d} = \{\gep, \gd-\gep\},\quad\quad \Delta^{+, \rm d}_1 = \{\gd, \gd\pm \gep\},\]
\[ \Pi^{\rm a} = \{\gd, \gep-\gd\},\quad\quad \Delta^{+, \rm a}_1 = \{\gd, \gep\pm \gd\}.\] Each of these corresponds to a Borel subalgebras 
$\fb^{{\rm b}}$ and a version $\gr^{{\rm b}}$ of $\gr$ (for $b=a,d$) as defined in \eqref{b3}. 
The latter are, by computation 
\[ \gr^{\rm d} = (\gep- \gd)/2,\quad\quad \gr^{\rm a}= (\gd-\gep)/2.\]
Note that $\gd-\gep$ is a root of $\fb^{\rm d}$, but instead  $\gep- \gd$ is a root of $\fb^{\rm a}$. Apart from this  $\fb^{\rm d}$ and  $\fb^{\rm a}$ have the same roots. Thus $\fb^{\rm d}$ and  $\fb^{\rm a}$ are related by  an odd reflection. We define the symmetric bilinear form $(\;,\;)$ on $\fh^*$ by $(\gep,\gep) = -(\gd,\gd) =1$ and $(\gd,\gep)=0$.
Up to a non-zero scalar multiple, this bilinear form agrees with the one defined in \eqref{b6} by  a general construction.  
 We denote the Verma module with highest weight $\gl$ induced from $\fb^{{\rm b}}$ by $M^{{\rm b}}(\gl)$. There are two versions of the dot action, which we denote by $(w, \gl) \lra w\cdot_{{\rm b}} \gl = w(\gl+\gr^{{\rm b}}) -\gr^{{\rm b}}$. The Borel subalgebra $\fb^{\rm d}$ is called {\it distinguished} in \cite{Kac2}.  We call   $\fb^{\rm a}$ {\it anti-distinguished}.
\section{A Family of  Elusive Cases.}\label{elc} Consider the isotropic root $\gc=\gep+\gd$, and for a positive integer $n,$ define $\gl_n , \mu_n \in \fh^*$ by $$\gl^{\rm d}_n + \gr^{\rm d}=(2n-1)\gc/2, \quad\mu^{\rm d}_n +\gr^{\rm d} =(2n-1)(\gd-\gep )/2.$$  We have
\be \label{s1}(\gl^{\rm d}_n+ \gr^{\rm d}, \gd^\vee) =(\gl^{\rm d}_n+ \gr^{\rm d},\gep^\vee) =(\gl^{\rm d}_n+ \gr^{\rm d}, \gep-\gd) =2n-1, \; (\gl^{\rm d}_n+ \gr^{\rm d}, \gc) =0.\ee
Thus 
$s_\gep(\gl^{\rm d}_n + \gr^{\rm d}) =\mu^{\rm d}_n +\gr^{\rm d}.$  
In this Section we fix $n$ and set $\gl^{\rm d} =\gl^{\rm d}_n$.   
Now if $v_{\gl^{\rm d}}$ is the highest weight vector in $M^{\rm d}(\gl^{\rm d})$, then   $v_{\gl^{\rm a}}:=e_{\gep-\gd}v_{\gl^{\rm d}}$ is a highest weight vector for $\fb^{\rm a}$ 
with weight 
$\gl^{\rm a}= \gl^{\rm d} + {\gep-\gd}$,
which generates a copy of the Verma module $M^{\rm a}(\gl^{\rm a})$.  Since 
$(\gl^{\rm d}+ \gr^{\rm d}, \gd-\gep) \neq 0$,  
 $e_{\gd-\gep}v_{\gl^{\rm a}}$ is a non-zero multiple of   $v_{\gl^{\rm d}}$,
and we have      
$U(\fg)v_{\gl^{\rm a}}= M^{\rm d}(\gl^{\rm d})$.  
Thus we may write 
\be \label{s7} M(\gl):=M^{\rm d}(\gl^{\rm d})=M^{\rm a}(\gl^{\rm a}).
\ee
 Moreover 
 we have $\gl^{\rm d}+ \gr^{\rm d}= \gl^{\rm a}+ \gr^{\rm a},$ compare \cite{M} Corollaries 8.6.2 and 8.6.3.
In this case $M(\gl)$ has length 8, and is multiplicity free.  For ease of comparison with \cite{Mas},
we use the same  notation for simple modules.  Thus the composition factors of 
$M(\gl)$ have the form $L_x$ for $x \in \{a, b, c,d, e, f, g,h \}$, and in the diagrams below we label a length one interval by $x$  if the corresponding 
subfactor is
isomorphic to   $L_x$. 
\\ \\
For a positive root $\eta,$ let $\cH_{\eta,r} =\{\nu\in \fh^*|(\nu+\gr,\eta)=r(\eta,\eta)/2 \}$. Let $r$ be a positive integer.  If $\eta$ is isotropic assume that $r=1$, and if $\eta$ is odd non-isotropic assume that $r$ is odd. 
A {\it \v Sapovalov element} $\gth_{\eta,r} \in U(\fb^-)$, for the pair $(\eta,r)$ 
has the property that  for $\nu\in \cH_{\eta,r}$, $\gth_{\eta,r} v_\nu \in M(\nu)$
is a non-zero highest weight vector.  A { \v Sapovalov element} $\gth_{\eta,r}$ is only determined modulo a left ideal  in $U(\fb^-)= U(\fn^-)\ot U(\fh)$. However for $\nu\in \cH_{\eta,r}$, the evaluation $\gth_{\eta,r}(\nu) \in U(\fn^-)$ is uniquely determined, and we have 
 $\gth_{\eta,r}(\nu)v_\nu= \gth_{\eta,r} v_\nu$. For a simple root $\eta$, we have
$\gth_{\eta,r} = e_{-\eta}^r$. 
If $\eta$ is isotropic set $ \gth_{\eta}=\gth_{\eta,1}$ and $ \cH_{\eta}=\cH_{\eta,1}$. 
\bl \label{k9} Let $\gc$ be  a positive isotropic root and suppose   
$\nu\in \cH_\gc$. Then \bi \itema There is a factor module $M^{\gc}(\nu)$ of $M(\nu)$ such that in $K(\cO)$ we have $[M(\nu)] =[M^{\gc}(\nu)] +[M^{\gc}(\nu - \gc)]$.
\itemb
If $K(\nu)$ is the kernel of the natural map 
$M(\nu)\lra M^{\gc}(\nu)$, then 
$[K(\nu)]=[M^{\gc}(\nu-\gc)].$
\itemc $\gth_{\gc}v_\nu \in K(\nu)$.
\ei \el 
\bpf (a) is shown in \cite{M21} Section 6 using \v Sapovalov elements, (b) is an immediate consequence and (c) is proved in \cite{M21} Corollary 6.8.
\epf \noi
For $\nu^{\rm b} \in
\mathfrak{h}^*$ define
\by \label{rue} A(\nu^{\rm b})_{0}  &=&  \{ \alpha \in {\Delta}^+_{0} | \ga/2 \mbox{ is not a root and } (\nu^{\rm b} + \gr^{\rm b},
\alpha^\vee) \in \mathbb{N} \backslash \{0\} \}\nn \\
A(\nu^{\rm b})_{1} &=& \{ \alpha \in \Delta^+_{1}  | \ga \mbox{ is non-isotropic  and } 
(\nu^{\rm b} + \gr^{\rm b}, \alpha^\vee ) \in 2\mathbb{N} + 1 \} \\
 A(\nu^{\rm b}) &=& A(\nu^{\rm b})_{0} \cup  A(\nu^{\rm b})_{1}.\nn \\ B(\nu^{\rm b}) &=& \{ \alpha \in {\Delta}^+_{1} | \ga \mbox{ is isotropic  and }  (\nu^{\rm b} + \gr^{\rm b},\alpha) = 0 \} .\nn 
\ey
\bl  \label{hon} Suppose  $\ga\in A(\nu^{\rm b})$ is a simple non-isotropic root of the Borel subalgebra $\fb^{{\rm b}}$, and $(\nu+\gr,\ga^\vee) =m$. Set $\mu^{\rm b} = s_\ga \cdot_{\rm b}\nu^{{\rm b}}$.  
Then \bi \itema  $e_{-\ga}^m v_{\nu^{\rm b}}$ generates a submodule of $M(\nu^{\rm b})$ which is isomorphic to $M^{{\rm b}}(\mu^{{\rm b}})$.
\itemb We have 
$\dim \Hom(M^{{\rm b}}(\mu^{{\rm b}}), M^{{\rm b}}(\nu^{{\rm b}}))=1$.\ei \el
\bpf Combine Lemma 9.2.1 and Theorem 9.3.2 from \cite{M}. \epf
\noi 
A key tool is the Jantzen sum formula, which is best expressed in the Grothendieck group $K(\cO)$ of the module category $\cO$. 
\ff{This is defined as the category of $\Z_2$-graded modules which belong to the usual BGG category $\cO$ when regarded as $\fg_0$-modules, see  \cite{H2}.  
}  
Objects in this category have finite  length, so we can define  $K(\cO)$
to be the free abelian group on the symbols  $[L]$, where $L$ is a simple module. If $M$ is any object of $\cO$ we define $[M]$ to be $\sum_L |M:L|[L]$
where $|M:L|$ is the multiplicity of $L$ as a composition factor of $M$.  
We note that 
$K(\cO)$ has a natural partial order.  For $A, B \in \cO$ we write $[A]
\ge [B]$ if $[A]-[B]$ is a linear combination of classes of simple modules with non-negative integer coefficients.  Clearly if $B$ is a subquotient of $A$
we have $[A] \ge [B]$.
\\ \\
Now in a  general Verma module $M(\nu)$ has a {\it Jantzen filtration} 
$M_1(\nu) \supseteq M_2(\nu) \supseteq  \ldots$ where  
$M_1(\nu)$ is the unique maximal submodule of  $M(\nu).$ The following is  the {\it Jantzen sum formula} (for the distinguished Borel subalgebra).
\bt \label{Jansum}
For all $\nu^{\rm d} \in \mathfrak{h}^*$
\be \label{lb} \sum_{i > 0} [{M}_{i}(\nu^{\rm d})] = \sum_{\alpha \in A(\nu^{\rm d})}[{M}(s_{\alpha}\cdot_{\rm d}\nu^{\rm d})] +
\sum_{\gc \in B(\nu^{\rm d})} [M^{\gc}(\nu^{\rm d} -\gc)].\ee
\et
\bpf See
\cite{M} Theorem 10.3.1 for a preliminary version and \cite{M21} Theorem 1.11 
for the statement in this form.  
\epf \noi
In the case of $\fg=\osp(3,2)$, set $\gc = \gep+\gd$. Then for                  $\gl=\gl^{\rm d}$ as above we have by \eqref{s1}, 
$A(\gl^{\rm d}) =\{\gep, \gd\}$ and $B(\gl^{\rm d}) =\{\gc\}$, so  \eqref{lb} 
 takes the form 
\by \label{hib} \sum_{i\ge 1}[M_i(\gl)]&=& [A] + [M^{\rm d}(\gs_\gd\cdot_{\rm d}\gl^{\rm d})] + [M^{\gc}(\gl-\gc)]\nn\\
&=& [L_b] +[L_c] +[L_e] +2[L_g] +2[L_d] +2[L_f] +3[L_h],\ey
see \cite{Mas} Equations (5.2.7) and   (5.2.8),  
where $A =M^{\rm d}(\gs_\gep\cdot_{\rm d}\gl^{\rm d})$ and $M(\gl)/M_1(\gl)\cong L_a:=L(\gl)$.  By Lemma \ref{hon}, 
$A$ is a join irreducible sumodule of $M(\gl)$,  but 
it is not clear that the other terms on the right of \eqref{hib} are even submodules. However 
we can improve the formula using Lemma \ref{k9} (b) and Corollary \ref{ao}. 
\bl \label{hz}The composition factors of each $M_i(\gl)$ are given as follows
\be \label{hb}M_4(\gl)=0,\quad  [M_3(\gl)]=[L_h],\quad [M_2(\gl)] =  [L_g] +[L_d] +[L_f] +[L_h]
\ee and
\be [M_1(\gl)] =[L_b] +[L_c] +[L_e] +[L_g] +[L_d] +[L_f] +[L_h].\ee \el
\bpf See \cite{Mas} Theorem 5.2.17.  We can also argue directly from   \eqref{hib} as follows.  Since no composition factor has multiplicity greater than 3, we have $M_4(\gl)=0$. Then 
$[L_h]$
is the only  composition factor with muliplicity 3 and so $[M_3(\gl)]=[L_h]$. More generally for any $j \ge 1$,  $[M_j(\gl)] = \sum_{x: \sum_{i\ge 1}[M_i(\gl): L_x] \ge j}[L_x]$ and this gives the result.  
\epf
\bl \label{ap} \bi \itema $[A] = [L_e] +[L_g]  +[L_f] +[L_h]$, 
\itemb $ [M^{\rm d}(\gs_\gd\cdot_{\rm d}\gl^{\rm d})] = [L_c]  +[L_g] +[L_d]  +[L_h]$, 
\itemc $ [M^{\gc}(\gl-\gc)] = [L_b]  +[L_d] +[L_f] +[L_h].$
\ei\el 
\bpf See \cite{Mas} page 73. \epf 
\bc \label{ao} $$[M^{\rm d}(\gs_\gd\cdot_{\rm d}\gl^{\rm d})]=[M^{\rm a}(\gs_\gd\cdot_{\rm a}\gl^{\rm a})].$$
\ec
\bpf  The composition factors of $M^{\rm d
}(\gs_\gd\cdot_{\rm d}\gl^{\rm d})$ are given by the Lemma. We use Corollary 4.8.6 with $\gl$ replaced by $\gs_\gd\cdot_{\rm a}\gl^{\rm a}$ from \cite{Mas} to find the composition factors of  $M^{\rm a}(\gs_\gd\cdot_{\rm a}\gl^{\rm a})$. The Corollary gives the highest weights of the composition  with respect to the Borel subalgebra $\fb^{\rm a}$. To complete the calculation we use \cite{Mas} Theorem 5.2.12 
 which gives the highest weight of each $L_x$ for both Borels.  Thus from the Corollary we know that the composition factors of $M^{\rm a}(\gs_\gd\cdot_{\rm a}\gl^{\rm a})$ are as follows:
\bi \itema $L^{\rm a}(\gs_\gd\cdot_{\rm a}\gl^{\rm a})\cong L_c$
\itemb $L^{\rm a}(\gs_\gd\cdot_{\rm a}\gl^{\rm a}-(\gep-\gd)) \cong L^{\rm a}(\gs_\gd\cdot_{\rm a}(\gl^{\rm a}-(\gep+\gd)) \cong L_d$
\itemc $L^{\rm a}(\gs_\gep\gs_\gd\cdot_{\rm a}\gl^{\rm a}-(\gep+\gd)) \cong  L^{\rm a}(\gs_\gep\gs_\gd\cdot_{\rm a}(\gl^{\rm a}+(\gep+\gd))) \cong  L_g$
\itemd $L^{\rm a}(\gs_\gep\gs_\gd\cdot_{\rm a}\gl^{\rm a})\cong L_h.$
\ei
\epf  
\noi
Now we state an improved version of the Jantzen sum formula.  Set $A =M^{\rm d}(\gs_\gep\cdot_{\rm d}\gl^{\rm d})$, $B =M^{\rm a}(\gs_\gd\cdot_{\rm a}\gl^{\rm a})$ and $C=K(\gl)$. Then $A, B, C$ are submodules of $M(\gl)$, with 
$A, B$ join irreducible, and we have 
 \by \label{heb1} \sum_{i\ge 1}[M_i(\gl)]&=& [A] + [B] + [C]\nn\\
&=& [L_b] +[L_c] +[L_e] +2[L_g] +2[L_d] +2[L_f] +3[L_h],\ey
where the $L_i$ are as before.
\\ \\
 By Lemma \ref{hz}, 
 $M(\gl)$ has simple socle $L_h$, and thus  $L_h$ is also the socle of  $A, B$ and $ C$. Now $A$ and $B$ are also join irreducible, with distinct simple factors, 
which are also different from $L_d, L_f$, $L_g$  and $L_h$. Thus we can choose the notation so that  $A$ and $B$ have unique simple factors $L_e$ and $L_c$ respectively.   This is consistent with Lemma \ref{ap}.
\bl \label{ly}The lattice of submodules of $A$ and $B$ have Hasse diagrams
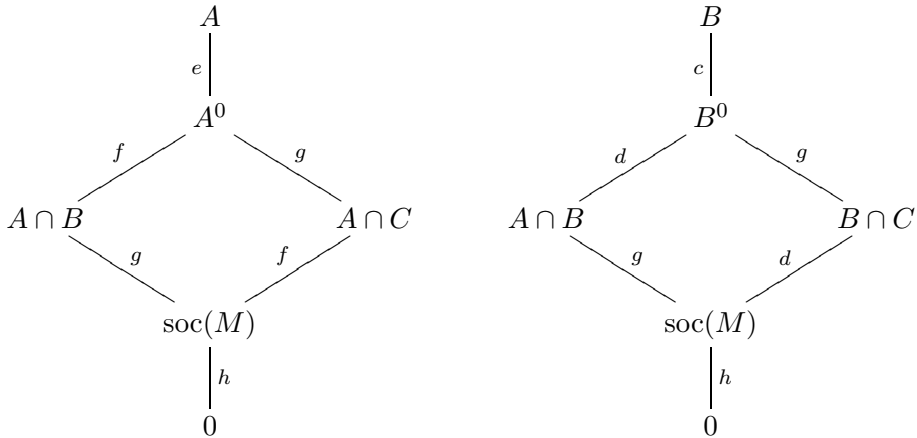
\Bc\[
\xymatrix{
&A\ar@{-}[d]_e&\\
&A^0&\\
A\cap B \ar@{-}[dr]^g \ar@{-}[ur]^f
 &&
A\cap C\ar@{-}[dl]_f \ar@{-}[ul]_g&\\
& \soc(M) \ar@{-}[d]^h & \\
&0&}
\xymatrix{
&B\ar@{-}[d]_c&\\
&B^0&\\
A\cap B \ar@{-}[dr]^g \ar@{-}[ur]^d
 &&
B\cap C\ar@{-}[dl]_d \ar@{-}[ul]_g&\\
& \soc(M) \ar@{-}[d]^h & \\
&0&}
\] \Ec
\el 
\bpf This follows from the following Lemma.\epf
\noi For $X$ a non-zero submodule of $M$ containing $\soc(M)$, set $\bar X= X/\soc(M)$.  
\bl \label{lo} We have \bi \itema $[A\cap B] = [L_g] +[L_h]$.  
\itemb $[A\cap C] = [L_f] +[L_h]$. 
\itemc $[B\cap C] = [L_d] +[L_h]$. 
\ei\el
\bpf The common composition factors of $A$ an $B$ are $L_g$ and $L_h$. If (a) is false, then 
$A\cap B = \soc(M)$.  But then $\overline{A+ B} = \bar A \op \bar B$ would have $L_g$ as a composition factor of multiplicity 2. The proofs of  (b), (c) are similar, since we know the composition factors of $C$. 
\epf \noi 
\noi Concerning  $L(C)$ we know that $C$ has simple  socle $\soc(M)$, and the composition factors of $C$ from Lemma \ref{ap}.  Thus  $L(\bar C)$ is isomorphic to one of the lattices from Lemma \ref{hag}. 
\bt \label{lp} The lattice of submodules  $L(\bar C)$ is isomorphic to the lattice $\gL$ in Case 2 or Case 4 of Lemma \ref{hag}. Thus the lattice of submodules  $L(C)$ is isomorphic to the lattices augmented $\gL$. \et
\bpf In Cases 1 and 5 $L(\bar C)$ has a unique atom, which by Lemma \ref{lo} would have to equal both $(\bar A \cap \bar C
)$ and $   (\bar B\cap \bar C) $. This would contradict the Lemma.
We still have to rule out Case 3, where $\gL$ is isomorphic to $\B$. We label 
$L(\bar C)$ with the join irreducibles we already know about.  
There is an extra  join irreducible, which we label $\bar D$.  
 From the information on composition factors in Lemmas \ref{ap} and \ref{lo}, we see that $\bar D \neq \bar A, \bar B, \bar C
$. If $\bar D=  \bar B\cap \bar C
,$ then $\bar C
=(\bar A\cap \bar B) +    (\bar A \cap \bar C
)   +   (\bar B\cap \bar C
) $ has composition factors 
$ L_d, L_f$ and  $L_h$. 
This contradicts Lemma \ref{ap}.
Thus $\bar D$ is a new join irreducible, and it follows that $M/\soc(M)$ has at least 8 composition factors, a contradiction. 
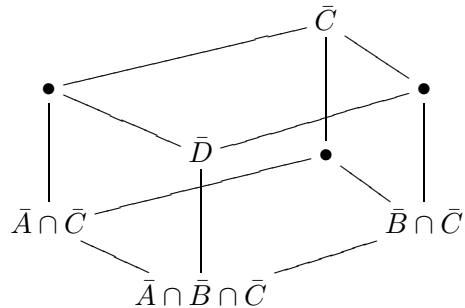
\[
\xymatrix@C=1pc@R=1pc{
& &  \bar C
& \\
\bullet  \ar@{-}[dr]
\ar@{-}[urr]
&&&
\bullet
\ar@{-}[dll]  
\ar@{-}[ul]
&\\
& \bar D&  \bullet
\ar@{-}[uu] \\
\bar A\cap  \bar C
 \ar@{-}[dr] \ar@{-}[urr]\ar@{-}[uu] &
&&
\bar B\cap \bar C
 \ar@{-}[dll] \ar@{-}[uu] \ar@{-}[ul]&\\
&  \bar A\cap \bar B \cap \bar C
 \ar@{-}[uu]
& \\
}
\] 
\epf 
\bl \label{y2}  \bi\itema If the Hasse diagram for  $L(\bar C)$ is as in Case 4 of Lemma \ref{hag}, then the  Hasse diagram  for $L(C)$  can be labelled as below on the left.
\itemb If the Hasse diagram for  $L(\bar C)$ is as in Case 2 of Lemma \ref{hag}, then there is a join irreducible submodule $D$ of $C$ such that the Hasse diagram  for $L(C)$  can be labelled as below on the right, possibly with the labels $A\cap C$ and 
$B\cap C$ interchanged. \ff{By Theorem \ref{kk}, in (b) we have $C = D +(B\cap C).$ Thus there is no need to interchange the labels, and the intervals are correctly labelled.}
\ei
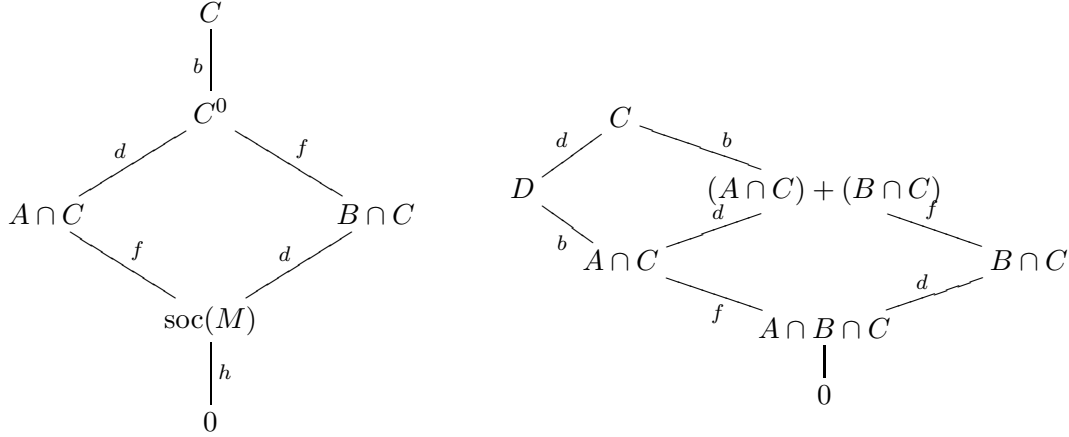
\[
\xymatrix{
&C\ar@{-}[d]_b&\\
&C^0&\\
A\cap C \ar@{-}[dr]^f \ar@{-}[ur]^d
 &&
B\cap C\ar@{-}[dl]_d \ar@{-}[ul]_f&\\
& \soc(M) \ar@{-}[d]^h & \\
&0&}
\xymatrix@C=1pc@R=1pc{
&&&\\
&&&\\
&C&\\
D\ar@{-}[dr]_b \ar@{-}[ur]^d &&
(A\cap C)+(B\cap C)\ar@{-}[dl]_d  
\ar@{-}[dr]^f\ar@{-}[ul]_b&\\
& A\cap C& &B\cap C \ar@{-}[dl]_d \\
 &&
A\cap B\cap C \ar@{-}[ul]^f&\\&&0\ar@{-}[u]&
} 
\] 
\el 
\bpf (a)  follows from Lemma \ref{lo}, If  (b) is false then $D$ would not be join irreducible, and it is easy to show that $C$ would not have simple socle.
\epf \noi 
\bc \label{lx}  There are only two ways the 
 Hasse diagrams for the lattice of submodules for 
$A,  B ,  C$ can fit together as a subdiagram of the Hasse diagram for  $N= M_1(\gl)/\soc(M)$.
\bi\itema If $(a)$ holds in Lemma \ref{y2}, then  the Hasse diagrams  fit together as in Step 1 of Subsection \ref{lz}.
\itemb If $(b)$ holds in Lemma \ref{y2}, then 
 the Hasse diagrams fit together as in Step 1 of Subsection \ref{lq}.
\ei
\ec
\bpf In each case, we know the Hasse diagrams for the lattice of submoudles of $A,  B ,  C$, and we know how these submodules intersect.
\epf \noi
Henceforth, unless otherwise stated, 
all  highest weights  are highest weghts for the distinguished Borel subalgebra $\fb$ and all
Verma modules are induced from $\fb$. We have $\gr= (\gep - \gd)/2$.
\noi Set 
$N_n= M_1(\gl_n)/\soc(M(\gl_n))$ and 
$C=K_n$ is the kernel of the natural map 
$M(\gl_n)\lra M^{\gc}(\gl_n)$.  The case $n=1$ is easily dealt with.
\bt \label{kk} If $n=1$, and 
$C=K_1$, then  
we have $C=D+(B\cap C)$ 
 in Lemma \ref{y2}, where $D = \gth_{\gc}v_{\gl_1}$ and 
the Hasse diagram for  $N_1$  has a subdiagram as in Step 1 of Subsection \ref{lq}.\et
\noi \bpf 
 By Lemma \ref{k9} $ D \subseteq K_{1}$. Note that $\gth_{\gc}v_{\gl_1}$ is a highest weight vector with weight ${\gl_1} -{\gc}= -{\gep}$.  Thus any weight of $D$ belongs to the set $-{\gep} -\sum_{\ga\in \Pi^{\rm d}}\N \ga $.  But by Lemma \ref{ly}, $(B\cap C)/\soc(M)\cong L_d,$ and by \cite{Mas} Theorem 5.2.12, 
$L_d = L(s_\gd\cdot_{{\rm d}}\gl_1)$. Now 
$\gl_1=
(\gep +\gd)/2-(\gep - \gd)/2 =\gd$, and $s_\gd\cdot_{{\rm d}}\gl_1 =0$. This is not a weight of $D$ by the above remark. Thus $L_d$ cannot be a composition factor of $D$, and it follows that $C= D+(B\cap C)$ is not join  irreducible as in (b) of Lemma \ref{y2}. The last statement follows from Corollary \ref{lx}.
\epf \noi
Since both cases actually occur in in Corollary \ref{lx}, this is as far as we can expect lattice theory and the results from \cite{Mas} alone to take us. 
\section{The Cases $n\ge 2$.} \label{osp1}
\bt \label{k7}For $n\ge2$,
$K_n$ is a highest weight submodule of $M(\gl_n)$ and 
the Hasse diagram for  $N_n$  has a subdiagram
 as in Step 1 of Subsection \ref{lz}.
\et\noi Now set  
$B^n =M(\gl_n)$, $A^n =M(\mu_n)$.  Note that $s_\gep\cdot \gl_n = \mu_n$.  Also by 
Lemma \ref{hon}, $M(\mu_n)$ is isomorphic to the submodule of $M(\gl_n)$, generated by $e_{-\gep}^{2n-1}v_{\gl_{n}}$. Define a map $\gs_{n}: A^n \lra B^n$  by $\gs_n(xv_{\mu_{n}}) = xe_{-\gep}^{2n-1}v_{\gl_{n}}$ 
and  set  $Q^n =B^n/\gs_n(A^n) $.  

\bl \label{k1} We have an exact sequence of complexes.

\be\label{ec}
\left.{\begin{array}{ccccccc}
 &0  &   &0 &&0\\
  &\lda & & \lda&& \lda\\
{\bf A}^\bullet:&A^1& \lra & A^2 &\lra&  A^3 &\lra\\
 &\lda & & \lda&& \lda\\
{\bf B}^\bullet:&B^1& \lra & B^2 &\lra&  B^3 &\lra\\
 &\lda & & \lda&& \lda\\
{\bf Q}^\bullet:&Q^1& \lra & Q^2 &\lra&  Q^3 &\lra\\
 &\lda & & \lda&& \lda\\
&0&  & 0&&0 
\end{array}}
 \right.
\ee
\el
\bpf The maps in the complexes $
{\bf A}^\bullet, 
{\bf B}^\bullet$ are given by  $$\gt_n: A^n  \lra A^{n+1}, \; xv_{\mu_{n}} \lra x  e_{\gep-\gd} v_{\mu_{n+1}}, \mbox{ and }  \psi_n: B^n  \lra B^{n+1}, \; xv_{\gl_{n}} \lra x  \gth_\gc v_{\gl_{n+1}}.$$
To check that the square 
$$
\left.{\begin{array}{ccc}A^n& \lra & A^{n+1} \\
 \lda &&  \lda\\
B^n& \lra & B^{n+1} \\
\end{array}}
 \right. 
$$ commutes, 
we show that  $
\psi_n \gs_{n}(v_{\mu_{n}}) =   \gs_{n+1}
\gt_n(v_{\mu_{n}})$. 
Now 
$\psi_n \gs_{n}(v_{\mu_{n}}) = e_{-\gep}^{2n-1} \gth_\gc v_{\gl_{n+1}}$ and  
$ \gs_{n+1}\gt_n(v_{\mu_{n}})= e_{\gep-\gd}e_{-\gep}^{2n+1}v_{\gl_{n+1}} $. 
These are equal by \cite{M21} (3.7). The map $\bar \psi_n: Q^n  \lra Q^{n+1}$ is induced by $\psi_n$. By  \cite{M21} Theorem 5.1, if $\eta$ is an isotropic root
and $\nu\in \cH_{\eta}$, then $ \gth_{\eta}^2 v_\nu =0$.  Thus the rows in the  diagram are complexes.  
\epf

\bl \label{k2} 
The complex ${\bf A}^\bullet$ is exact at $A^n$ for all $n\ge2$.
\el
\bpf Since ${\gep-\gd} $ is the negative of a simple root of $\fb$ this follows easily from the PBW theorem.\epf\noi 
Now to study the complex ${\bf Q}^\bullet$, 
define 
$$
\Gd^\fq =\{ -\gep, -\gep + \gd, \gd, 2\gd, \gep + \gd, \gep\}, \quad \Gd^\fm =\{ \gep - \gd, -\gd, -2\gd, -\gep - \gd\},$$ 
$$\fq = \fh \op\bigoplus_{\eta\in \Gd^\fq}\fg^\eta, 
\quad\fm=\bigoplus_{\eta\in\Gd^\fm}\fg^\eta .
$$  Then $\fq$ is a parabolic subalgebra of $\fg$, whose Levi part $\fl$ satisfies $$\fk:=[\fl,\fl] =\span \{h_{\gep}, e_{\gep}, e_{-\gep}\}.$$ 
Also $\fg= \fm\op\fq$ and $\fq= \ttk {e_{-\gep}}\op\fb$.  Now, since  $(\gl_n ,\gep) = n-1$, $U(\fq)\ot_{U(\fb)}\ttk v_{\gl_n}$ maps onto the simple 
$U(\fq)$-module, $L_n$ with character $\tte^{n\gd}(\tte^{n\gep} - \tte^{-n\gep} )/(\tte^{\gep} - \tte^{-\gep} )$. In particular $L_n$ has highest weight 
$\gl_n= (n-1)\gep +n\gd$ with highest weight vector $w_{\gl_n}.$ 
\bl \label{k3} $Q^n \cong \Ind_\fq^\fg \; L_n$.
\el \bpf \cite{Mas} Lemma 5.2.3. For convenience we sketch the proof.
First there is a surjection  from $B^n=M(\gl_n) $ onto $\Ind_\fq^\fg \; L_n$  sending $v_{\gl_n}$ to $w_{\gl_n}$. But the image of $A^n$ in $B^n$ is generated by $e_{-\gep}^{2n-1}v_{\gl_{n}}$ which maps to $e_{-\gep}^{2n-1}w_{\gl_{n}} = 0$ in $\Ind_\fq^\fg \; L_n$.  Thus 
 $Q^n\cong B^n/A^n$ maps onto $\Ind_\fq^\fg \; L_n$. On the other hand $B^n, A^n$ and $\Ind_\fq^\fg \; L_n$ are all induced from modules with known characters, so it is a simple matter to compute their characters, and show that $Q^n$ and $\Ind_\fq^\fg \; L_n$ have the same characters.
\epf \noi  

\bl \label{h4} Let $\ft =\span \{h_{\gd},  e_{\pm\gd}, e_{\pm2\gd} \} \cong \osp(1,2)$. 
\bi \itema
As a $U(\ft)$-module, $Q^n$ has a Verma flag.  \itemb Any non-zero submodule of $Q^n$ has GK\ff{Gelfand-Kirillov, see \cite{KrLe}}-dimension one.\ei\el
\bpf Since $\fg= \fm\op\fq$, $Q^n =U(\fm)\ot L_n$ is a free $U(\fm)$-module. Set $\fr = \span \{ e_{-\gd}, e_{-2\gd} \}$.  Then $\fr \subset \fm$, and $U(\fm) = \op_{i=1}^4 U(\fr)z_i$ is a free $U(\fr)$-module of rank 4, where 
$$z_1= 1,\quad z_2=e_{\gep-\gd},\quad  z_3= e_{-\gd-\gep}  \mbox{ and } z_4= e_{-\gd-\gep} e_{\gep-\gd}.$$ On the other hand $L_n= \span \{ e_{-\gep}^j w_{\gl_n} \}_{j=0}^{2n-2}$, and it follows that 
\be\label{h2} Q^n= \op_{i=1}^4\op_{j=0}^{2n-2}U(\fr)z_i e_{-\gep}^{j}
w_{\gl_n} .\ee Now define for $k\ge 1$, 
\[T_k =\{ (i,j)| e_\gd^kz_i e_{-\gep}^{j}
w_{\gl_n}=0  \} \mbox{ and } M_k = \bigoplus_{(i,j) \in T_k}U(\fr)z_i e_{-\gep}^{j}
w_{\gl_n}.\] Since $e_\gd$ acts nilpotently on $Q^n$, there is a least integer $m$  such that $M_m =Q^n.$ Thus we have a filtration 
$$0=M_0 \subset M_1 \subset \ldots \subset M_m=Q^n.$$
By construction  for $1 \le k \le m$, $M_k/M_{k-1}$ is a direct sum of $U(\fk)$-Verma modules, and this gives the result.\\ \\
(b) It  is well-known that  any non-zero submodule of  a $U(\fk)$-Verma module
GK-dimension one and the result follows from (a) and an easy argument, see \cite{M4}.
\epf

\bl \label{k4} For $n\ge 2$, the lattice of submodules of $Q^n$ has the form
\[
\xymatrix@C=0.8pc@R=0.8pc
{
&Q^n =U(\fg)w_{\gl_{n}} \ar@{-}[d]&\\
&V_1+V_2&\\
V_{1} \ar@{-}[dr] \ar@{-}[ur] &&
V_{2} \ar@{-}[dl] \ar@{-}[ul]&\\
& V_3= V_1\cap V_2\ar@{-}[d] & \\
&0&}
\]
with $V_i = U(\fg)v_i$ where 
$$v_1= e_{-\gd}^{2n-1}e_{\gep-\gd}w_{\gl_{n}},\quad v_2= \gth_\gc w_{\gl_{n}} \mbox{ and } v_3= e_{\gd-\gep} e_{-\gd}^{2n-1}e_{\gep-\gd}w_{\gl_{n}}.$$
In particular the highest weights of the composition factors of $Q^n$ 
are 
\be \label{h1} \gl_n,  \quad \gs_\gd  \cdot\gl_n -(\gep-\gd) , \quad \gl_n -\gc  \mbox{ and } \quad  \gs_\gd  \cdot \gl_n.\ee
\el
\bpf First we show that the $v_i$, $i=1,2,3$ generate proper submodules of $Q^n$.  For $v_2$ this follows from the general remarks about \v Sapovalov elements and the observation that  $\gth_\gc  v_{\gl_{n}}$ cannot belong to the submodule of $M({\gl_{n}})$ generated by $e_{-\gep}^{2n-1}v_{\gl_{n}}$ by weight considerations. Next note that $ e_{\gep-\gd}w_{\gl_{n}}$ is a highest weight vector for $\fb^{\rm a}$ which also generates $Q^n$, compare \eqref{s7}. Since $\gd$ is a simple root of $\fb^{\rm a}$, $v_1$ is a highest weight vector for  $\fb^{\rm a}$ compare Lemma \ref{hon}.  It is a singular vector for $\fb^{\rm d}$.   Then changing Borels again, we see that 
$v_3= e_{\gd-\gep} v_1$ is a highest weight vector for $\fb^{\rm d}.$ Note that the weights of $v_1 , v_2 , v_3$ are the last 3 weights listed in \eqref{h1}.  Because  the quadratic Casimir element $\gO$ acts on $Q^n$ as the scalar 
$(\gl_{n}+2\gr, \gl_{n})$, \cite{M} Lemma 8.5.3, it follows that if 
$\gl_{n} -\gz$ is the highest weight of a composition factor of  $Q^n$
we have $(\gz,\gz) =2(\gl_{n}+\gr,\gz)$, and $\gz$ is a weight of $U(\fm)$. 
Using these facts, it is not hard to show by a direct computation that the only possible highest weights of the composition factors are those listed in \eqref{h1}, see \cite{Mas}  Theorem 4.31. We remark that by \eqref{h2}, $Q^n$ is a finitely generated free module over the polynomial ring $\ttk[e_{-\gd}]$, so there are not many cases to check. In addition $Q^n$ is multiplicity free, since $B^n$ is multiplicity free by \eqref{hib} and Lemma \ref{hz}. 
 The only thing left to show is that $V_3 \subset V_2$.  If this is not the case then $V_2$ is simple.
The highest weight of $V_2$ with respect to $\fb^{\rm d}$
is $\nu^d = \gl_n -\gc$ and 
$$\nu^d+ \gr^d = (2n-3) \gc/2.$$  Since $(\gc,\gep-\gd)\neq 0$, $
e_{\gep-\gd}v_2$ is a highest weight vector for  $\fb^{\rm a}$ with weight $\nu^{\rm a}$ which also generates $V_2.$
By \cite{M}  Corollary 8.6.3, 
$\nu^{\rm d}+ \gr^{\rm d}=\nu^{\rm a}+ \gr^{\rm a}$, so 
$(\nu^{\rm a}+ \gr^{\rm a}, \gd^\vee) =
2n-3 $ is odd and positive.  Since $e_\gd$ is a simple root for  $\fb^{\rm a}$, and $V_2$ is a submodule of a free $\ttk[e_{-\gd}]$-module,  by Lemma \ref{h4} (and thus also free), it follows that  
$e_{-\gd}^{2n-3}e_{\gep-\gd}v_{2}$ is a non-zero highest weight vector in $V_2$.
Thus $V_2$ cannot be simple. In fact by a long and tedious computation 
$e_{-\gd}^{2n-3}e_{\gep-\gd}v_{2}$ is a non-zero scalar multiple of $v_3$. 
\epf
\bc \label{h5}  
 For $n\ge 1$,the maps $\bar \psi_n: Q^n  \lra Q^{n+1}$ 
in the complex $
{\bf Q}^\bullet$ are non-zero. 
\ec
\bpf This follows since  $\psi_n(v_{\gl_{n}})=\gth_\gc v_{\gl_{n+1}}$ and by the proof of the Lemma, this cannot belong to the submodule of $M({\gl_{n+1}})$ generated by $e_{-\gep}^{2n+1}v_{\gl_{n+1}}$.
\epf

\bl \label{k5} 
The complex 
${\bf Q}^\bullet$ is exact at $Q^n$ for all $n\ge2$.
\el 
\bpf Consider the diagram from Lemma  \ref{k4}. Clearly $\Im \bar \psi_{n-1} = V_2$. 
Since ${\bf Q}^\bullet$ is a complex, $\Im \bar \psi_{n-1}  \subseteq \Ker \bar \psi_{n}$.  If the inclusion is strict, then by Corollary \ref{h5} and Lemma \ref{k4}, 
$Q^{n}/\Ker \bar \psi_{n}$ is  the unique simple factor of $Q^{n}$, which has finite dimension.  
But then $\Im \bar \psi_{n}$
is  a finite dimensional submodule of $Q^{n+1}$
contradicting Lemma \ref{h4}. \epf
\bc \label{k6} 
The complex 
${\bf B}^\bullet$ is exact at $B^n$ for all $n\ge2$.
\ec
\bpf For $X= B$ or $Q$, denote the homology of the complex 
${\bf X}^\bullet$ is exact at $X^n$ by  $H(X^n)$. From the long exact homology sequence and Lemma \ref{k2}, we see that $H(B^n)=H(Q^n).$ Thus the result follows from Lemma  \ref{k5}.
\epf
\noi {\it Proof of Theorem \ref{k7}}. First suppose $n\ge1$.  In the complex $
{\bf B}^\bullet$ from \eqref{ec} we have $B^n =M(\gl_n)$, 
and the map $\psi_n:M(\gl_n)\lra M(\gl_{n+1})$
is given by $\psi_n(xv_{\gl_{n}})=  x\gth_\gc v_{\gl_{n+1}}$. Thus $\Im \psi_n= U(\fg) \gth_\gc v_{\gl_{n+1}}$.  Also by Lemma \ref{k9}, we have $ \Im \psi_{n-1} \subseteq K_{n}.$ Now if $n\ge2$, then by Corollary  \ref{k6} we have an exact sequence
$$0 \lra  \Im \psi_{n-1}
\lra
M(\gl_n)\lra \Im \psi_{n}\lra 0.$$ Thus in $K(\cO)$ we have, using Lemma \ref{k9}
$$
[M(\gl_n)]=[\Im \psi_{n-1}]
+[\Im \psi_{n}]\le  [K_{n}]+
[K_{n+1}] =
[M(\gl_n)].$$ Hence  $K_{n}=\Im \psi_{n-1} $ is an image of a Verma module, and this gives the result. 
\hfill  $\Box$

\section{Application: Building the Lattice of Submodules of the Verma  Module $M(\gl_n)$.} \label{uv2}
The set of submodules of a module $N$ becomes a lattice $L(N)$ when $\vee$ is taken to mean $+$, and $ \wedge$ to mean $\cap$. We want to  construct the lattice of submodules of  the Verma module $M(\gl_n)$ as in \eqref{s7}. Since $M(\gl_n)$  has a unique maximal submodule  $M_1(\gl_n)$ and a simple socle $M_3(\gl_n)$, an equivalent problem is 
to determine the lattice of submodules  of $N_n=M_1(\gl_n)/M_3(\gl_n)$. There are two possibilities, depending on which case holds in Corollary \ref{lx}.

\subsection{The Case $n=1$.}\label{lq} Suppose that $n=1$ and set 
$N=N_1$. We construct the Hasse diagram for $L(N)$, and show that 
  $L(N)$ has 21 elements. Thus $M(\gl_1)$ has 23 submodules. There are several unusual features  in this case. The submodules $A,B, C$ generate $L(N)$ but $C$ is not join irreducible.  The maximal submodule $D$ of $C$ is join irreducible, but the sublattice of $L(N)$ generated by  $A,B, D$ does not contain  $B\cap C$. The submodules $A, B,  D$ and  $B\cap C$ form a minimal generating set of join irreducibles. Thus  $L(N)$ has rank 4.  All this is apparent from Step1.
 \\ \\
{\bf Step 1.}
From Theorem \ref{kk}, we know that $N$ 
 contains at least the following submodules, where $A, B, D$ are join irreducible.  Furthermore any composition factor of $N$ is a also composition factor of $A, B$ or $ C$.
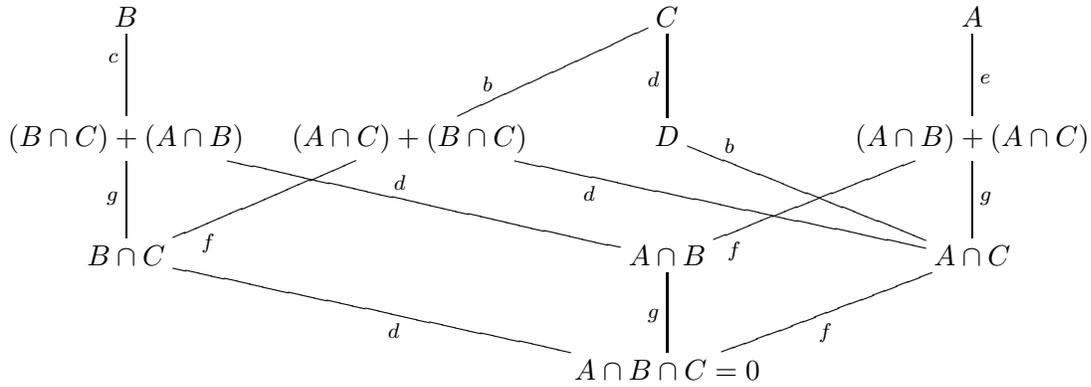
\[
\xymatrix@C=1pc@R=1pc{
 \\B
\ar@{-}[dd]_<<<<c&&C
&
&A
&&&&&&&\\ \\
 (B\cap C)+ (A \cap B)   \ar@{-}[ddrr]^d    
 &(A\cap C) +(B\cap C)\ar@{-}[uur]^<<<<<<{b} 
\ar@{-}[ddrrr]_<<<<<<<<<<<d
 &D\ar@{-}[uu]^d
\ar@{-}[ddrr]^<<<<<<b&& (A \cap B)+    (A\cap C)
\ar@{-}[uu]_e 
\\&&  &&&\\
B\cap C
 \ar@{-}[ddrr]_d \ar@{-}[uur]_<<<<<{f}\ar@{-}[uu]^g &
& A \cap B\ar@{-}[uurr]_<<<f &&A\cap C \ar@{-}[ddll]^f 
\ar@{-}[uu]_g&
\\&&&&&\\
&&  A\cap  B\cap C =0 \ar@{-}[uu]^{g}
& \\
}
\]
The diagrams can be drawn more compactly using the abbreviations
$$G=(B\cap C)+ (A \cap B), \quad F=(A \cap B)+    (A\cap C)$$
$$E=(A\cap C) +(B\cap C), \quad  S= (A \cap B)   +   (A\cap C)+(B\cap C). $$
Thus the Hasse diagram for $L(N)$ contains the following subdiagram.
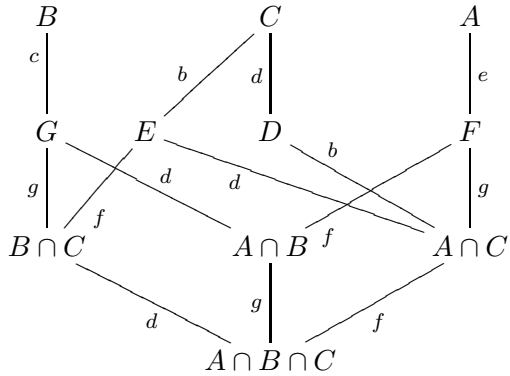
\[
\xymatrix@C=1pc@R=1pc{
 B
\ar@{-}[dd]_<<<<c&&C
&
&A
&&&&&&&\\\\
 G   \ar@{-}[ddrr]^d    
 &E\ar@{-}[uur]^<<<<<<{b} 
\ar@{-}[ddrrr]_<<<<<<<<<<<d
 &D\ar@{-}[uu]^d
\ar@{-}[ddrr]^<<<<<<b&& F
\ar@{-}[uu]_e 
\\
&&  &&&
\\
B\cap C
 \ar@{-}[ddrr]_d \ar@{-}[uur]_<<<<<{f}\ar@{-}[uu]^g
 &
& A \cap B\ar@{-}[uurr]_<<<f &&A\cap C \ar@{-}[ddll]^f 
\ar@{-}[uu]_g&
\\&&&&&\\
&&  A \cap B\cap C \ar@{-}[uu]^{g}
& \\
}\]
{\bf Step 2.} In Step 1, we have a box with one corner missing. 
Complete the box by adding $S=E+F=E+G=F+G$.  Also add  $D+F$. %
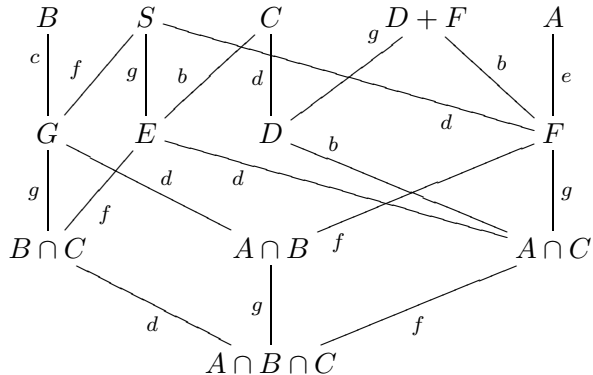
\[
\xymatrix@C=1pc@R=1pc{
B
\ar@{-}[dd]_<<<<c&S \ar@{-}[ddl]_>>>>>>{f} 
\ar@{-}[dd]_<<<<<<{g} 
\ar@{-}[ddrrr]_>>>>>>>>>>>>d&C
& D+F \ar@{-}[ddl]_<<<<g \ar@{-}[ddr]^b 
&A
&&&&&&&\\
\\
 G   \ar@{-}[ddrr]^d    
 &E
\ar@{-}[uur]^<<<<<<{b} 
\ar@{-}[ddrrr]_<<<<<<<<<<<d
 &D\ar@{-}[uu]^d
\ar@{-}[ddrr]^<<<<<<b&& F
\ar@{-}[uu]_e 
\\
&&&&&\\
B\cap C
 \ar@{-}[ddrr]_d \ar@{-}[uur]_<<<<<<{f}\ar@{-}[uu]^g &
& A \cap B\ar@{-}[uurr]_<<<f &&A\cap C \ar@{-}[ddll]^f 
\ar@{-}[uu]_g&
\\
&&&&&\\
&&  A \cap B\cap C \ar@{-}[uu]^{g}
& \\
}
\] 
\brs {\rm \bi
\itema In Step 2,  all modules in the row starting with $A$ have length 3.  
This will help with the rest of the construction.  We add some notes to show that  we have not missed any sums of submodules from previous steps.  At the end  we list the composition factors of  all non-simple submodules,  to ensure that nothing has been repeated.
\itemb It is easy to see how to proceed based on a consideration of composition factors.  For example consider the quadrilateral  in Step 3, with base $F$.  This submodule has composition factors $L_f$ and $L_g$, and is covered by two submodules $A$ and $S$, that also have $L_e$ or $L_d$ as a composition factor. This can be seen already in Step 2. By completing the quadrilateral in Step 3, we obtain the length 4 submodule $A+S$.  The length 4 submodule $A+D$ is constructed simililarly.    Opposite sides of a quadrilateral have the same  labels by the second isomorphism theorem.
  \ei}
\ers \noi
{\bf Step 3.}
Next add all submodules of length 4.
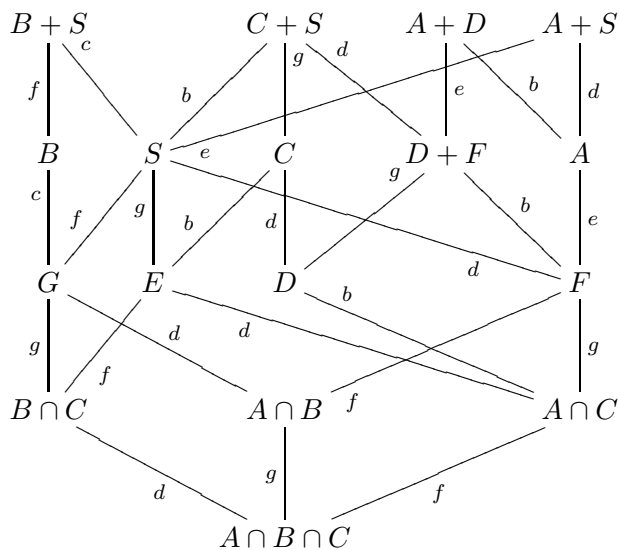
\[
\xymatrix@C=1.2pc@R=1.2pc{
B+   S\ar@{-}[dd]_{f}\ar@{-}[ddr]^<<<c 
 & & C+S  \ar@{-}[ddl]_>>>>>>b  \ar@{-}[dd]^<<<{g} 
\ar@{-}[ddr]^<<<<<d 
&A+D\ar@{-}[dd]^e
\ar@{-}[ddr]^>>>>>>>>b 
&A+S\ar@{-}[dd]^d \ar@{-}[ddlll]^>>>>>e
\\
&&
\\B
\ar@{-}[dd]_<<<<c&S \ar@{-}[ddl]_>>>>>>{f} 
\ar@{-}[dd]_<<<<<<{g} 
\ar@{-}[ddrrr]_>>>>>>>>>>>>d&C
& D+F \ar@{-}[ddl]_<<<<g \ar@{-}[ddr]^b 
&A
&&&&&&&\\
\\
 G   \ar@{-}[ddrr]^d    
 &E
\ar@{-}[uur]^<<<<<<{b} 
\ar@{-}[ddrrr]_<<<<<<<<<<<d
 &D\ar@{-}[uu]^d
\ar@{-}[ddrr]^<<<<<<b&& F
\ar@{-}[uu]_e 
\\
&&  &&&\\
B\cap C
 \ar@{-}[ddrr]_d \ar@{-}[uur]_<<<<<<{f}\ar@{-}[uu]^g &
& A \cap B\ar@{-}[uurr]_<<<f &&A\cap C \ar@{-}[ddll]^f 
\ar@{-}[uu]_g&
\\
&&&&&\\
&&  A \cap B\cap C \ar@{-}[uu]^{g}
& \\
}
\] 
\bi \itema The top row contains sums of 4 of the 5  length 3 submodules.  However  \itemb $F\subset A$, so $A+D+F = A+D$.  We give this submodule a shorter name.
\itemc $F \subset S$, so 
$C+F \subseteq C+S$, but $C+F \supseteq S$, so $C+F =C+S = C+D+F$.  
\itemd $C+S = D+F+S$
\iteme  $B+C= B+D+F,  A+C$  and $A+B$ have length 5.  
\itemf Together (a), (c), (d), (f) account for all 10 sums of the length 3 submodules.
\ei
{\bf Step 4.} 
Next add all submodules of length 5.
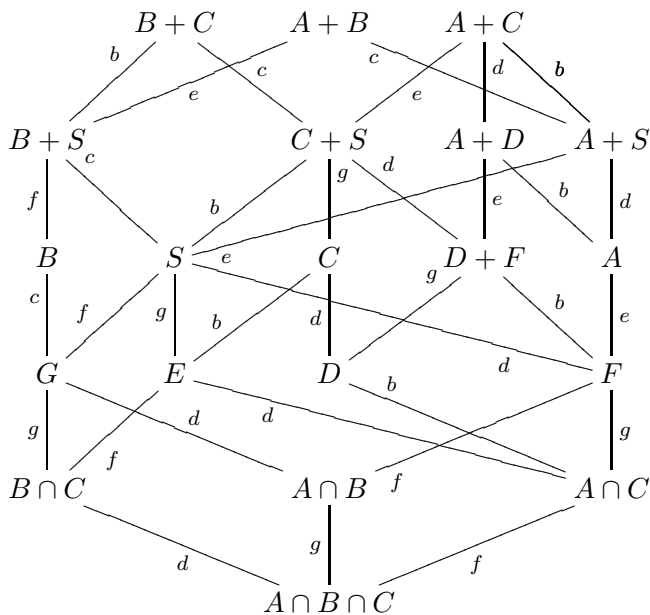
\[
\xymatrix@C=1pc@R=1pc{
&B+C\ar@{-}[ddl]_<<<<<<b \ar@{-}[ddr]^c 
&A+B
\ar@{-}[ddrr]_<<c \ar@{-}[ddll]^e
&
A+C\ar@{-}[dd]^<<<<d
\ar@{-}[ddr]^b \ar@{-}[ddr]^b\ar@{-}[ddl]^e
\\
&&&&&&&\\
B+   S\ar@{-}[dd]_{f}
\ar@{-}[ddr]^<<<c 
 & 
& C+S  \ar@{-}[ddl]_>>>>>>b  \ar@{-}[dd]^<<<{g} 
\ar@{-}[ddr]^<<<<<d 
&A+D\ar@{-}[dd]^e
\ar@{-}[ddr]^>>>>>>>>b
&A+S\ar@{-}[dd]^d \ar@{-}[ddlll]^>>>>>e
\\
&&
\\B
\ar@{-}[dd]_<<<<c&S \ar@{-}[ddl]_>>>>>>{f} 
\ar@{-}[dd]_<<<<<<{g} 
\ar@{-}[ddrrr]_>>>>>>>>>>>>d&C
& D+F \ar@{-}[ddl]_<<<<g \ar@{-}[ddr]^b 
&A
&&&&&&&\\ \\ G   \ar@{-}[ddrr]^d    
 &E
\ar@{-}[uur]^<<<<<<{b} 
\ar@{-}[ddrrr]_<<<<<<<<<<<d
 &D\ar@{-}[uu]^d
\ar@{-}[ddrr]^<<<<<<b&& F
\ar@{-}[uu]_e   
\\
& 
&  
&&&
\\
B\cap C
 \ar@{-}[ddrr]_d \ar@{-}[uur]_<<<<<<{f}\ar@{-}[uu]^g &
& A \cap B\ar@{-}[uurr]_<<<f &&A\cap C \ar@{-}[ddll]^f 
\ar@{-}[uu]_g&
\\
&&&&&\\
&&  A \cap B\cap C \ar@{-}[uu]^{g}
& \\
}
\]  \bi 
\itema The top row contains sums of 3 of the 4  length 4 submodules.  However
\itemb $B+C \supset S= (A \cap B)   +   (A\cap C)+(B\cap C).$ So $B+C+S= B+C.$ We give this submodule the shorter name.
\itemc
 We  do likewise with  $A+C=A+D+S$.
\itemd $A+C=(A+D)+(C+S) =(A+S)+(C+S) $.
\iteme $(B+S)+(A+S) =N$, since it contains all composition factors of $N.$
\itemf Together (a), (d), (e) account for all 6 sums of the length 4 submodules.
\ei
{\bf Final Step.}  Add the sum of the maximal submodules $A+  B +  C$
from Step 4.
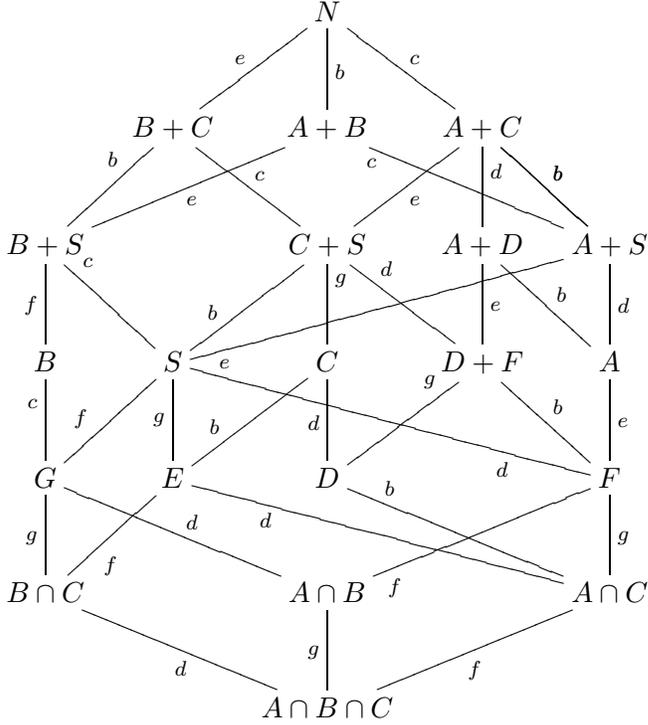
\[
\xymatrix@C=1pc@R=1pc{
 &&N
\ar@{-}[ddr]^c
\ar@{-}[ddl]_e
\ar@{-}[dd]^b
 &
\\&&&&&&&\\
&B+C\ar@{-}[ddl]_<<<<<<b \ar@{-}[ddr]^c 
&A+B
\ar@{-}[ddrr]_<<c \ar@{-}[ddll]^e
&A+C\ar@{-}[dd]^<<<<d 
\ar@{-}[ddr]^b \ar@{-}[ddr]^b\ar@{-}[ddl]^e
\\
&&&&&&&\\
B+   S\ar@{-}[dd]_{f}
\ar@{-}[ddr]^<<<c 
 & 
& C+S  \ar@{-}[ddl]_>>>>>>b  \ar@{-}[dd]^<<<{g} 
\ar@{-}[ddr]^<<<<<d 
&A+D\ar@{-}[dd]^e
\ar@{-}[ddr]^>>>>>>>>b
&A+S\ar@{-}[dd]^d \ar@{-}[ddlll]^>>>>>e
\\
&&
\\B
\ar@{-}[dd]_<<<<c&S \ar@{-}[ddl]_>>>>>>{f} 
\ar@{-}[dd]_<<<<<<{g} 
\ar@{-}[ddrrr]_>>>>>>>>>>>>d&C
& D+F \ar@{-}[ddl]_<<<<g \ar@{-}[ddr]^b 
&A
&&&&&&&\\
\\
 G   \ar@{-}[ddrr]^d    
 &E
\ar@{-}[uur]^<<<<<<{b} 
\ar@{-}[ddrrr]_<<<<<<<<<<<d
 &D\ar@{-}[uu]^d
\ar@{-}[ddrr]^<<<<<<b&& F
\ar@{-}[uu]_e   
\\
& 
&  
&&&
\\
B\cap C
 \ar@{-}[ddrr]_d \ar@{-}[uur]_<<<<<<{f}\ar@{-}[uu]^g &
& A \cap B\ar@{-}[uurr]_<<<f &&A\cap C \ar@{-}[ddll]^f 
\ar@{-}[uu]_g&
\\
&&&&&\\
&&  A \cap B\cap C \ar@{-}[uu]^{g}
& \\
}
\] 
\noi 
{\bf Composition Factors}.  To show that all submodules in the diagram are distinct we list the composition factors of non-simple submodules.  The equations below hold in the Grothendieck group $K(\cO),$ but we use a more compact notation than before.  Thus if $X$ is a submodule of $N$, we write $X= \sum_{|X:L_x| =1} x.$ Modules in the same row have the same length. 
$$G = d+g,\quad E= d+f,\quad D=b+f,\quad F=f+g $$ 
$$B= c+ d+g,\quad S=d+f+g,\quad C= b+d+f,\quad D+F = b+f+g, \quad  A=e+f+g.$$
$$B+ S =c+d+f+g, C+S= b+d+f+g, A+D=b+e+f+g, A+S=d+e+f+g,$$ $$B+C=b+d+c+f+g, A+ C=b+d+e+f+g, A+B= c+ d+e+f+g.$$

\noi
The submodule $S$ is the socle of $N$, and the Hasse diagram for $L(S)$ is the same as that of the box  lattice $\B$. To determine the the Hasse diagram for $L(N/S)$  consider all submodules not contained in $S$, and 
identify two submodules that are joined by en edge labelled by a submodule of $S$.  We find that  $L(N/S)\cong \B$ as lattices. 

\subsection{The Case $n\ge 2$.}\label{lz} Suppose that $n\ge 2$. We show that 
  $L(N_n)$ is isomorphic to the  restricted free distributive lattice of rank 3. Thus $L(M(\gl_n))$ is isomorphic to extended $\gL_3$. 
\\ \\
{\bf Step 1.}
From Theorem \ref{k7}, we know that $N=N_n$ contains at least the following submodules, where $A, B, C$ are join irreducible 
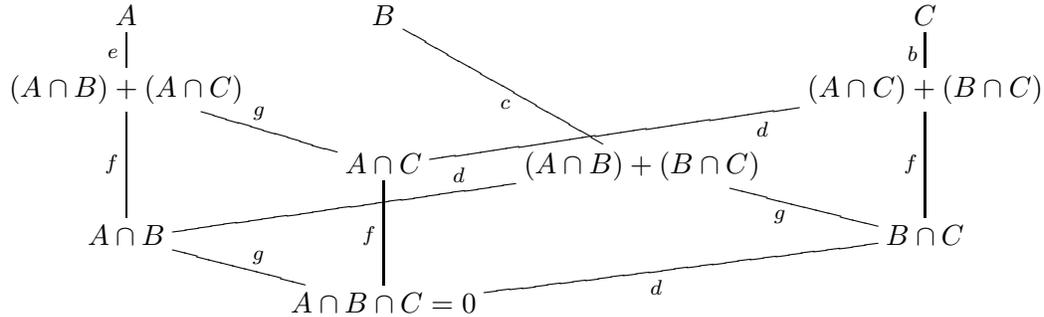
\[
\xymatrix@C=1pc@R=1pc{
A &B &  &C \\
 (A\cap B)+ (A \cap C)   \ar@{-}[dr]^g \ar@{-}[u]^e\ar@{-} &&&
  (A \cap C)+    (B\cap C)  \ar@{-}[u]^b\ar@{-}[dll]^<<<<<<{d} &\\
& A \cap C &  (A\cap B) +  (B\cap C)  \ar@{-}[uul]^c\\
A\cap B \ar@{-}[dr]^g \ar@{-}[urr]^>>>>>>>>{d}\ar@{-}[uu]^f  &
&&
B\cap C \ar@{-}[dll]^d \ar@{-}[uu]^f \ar@{-}[ul]^g&\\
&  A\cap B \cap C =0\ar@{-}[uu]^{f}
& \\
}
\] 
We  omit the remaining details, because as remarked in the first version of this paper, the result is predictable in advance, since from Step 1,  we see that the poset of join irreducble subsets of $M(\gl_n))$ is isomorphic to to $\B$, and $J(\B)$ is isomorphic to extended $\gL_3$. For full details on the stepwise  construction  of $L(N_n)$ see \cite{M3}.

\end{document}